\documentclass[a4paper,12pt]{amsart}
\usepackage{amsopn,amsfonts,amssymb,amsthm,mathtools,amsmath}
\usepackage{color}
\usepackage{latexsym}

\usepackage{wrapfig}
\usepackage[T1]{fontenc}
\usepackage{lmodern}
\usepackage{textcomp}
\usepackage[utf8]{inputenc}
\usepackage{enumerate}
\usepackage[shortlabels]{enumitem}
\setlist[enumerate, 1]{\sc(1)}
\usepackage{float}
\usepackage{hhline}
\usepackage{multirow}
\usepackage{anysize}
\usepackage{tikz}
\marginsize{2cm}{2cm}{2cm}{2cm}
\usepackage{graphicx}
\usepackage{circledsteps}
\usepackage{hyperref}
\hypersetup{colorlinks=true,linkcolor=blue,citecolor=blue}

\usepackage{enumerate}
\numberwithin{equation}{section}

\def\z{\mathfrak{z}}

\def\g{\mathfrak{g}}
\def\h{\mathfrak{h}}

\def\sl{\mathfrak{sl}}
\def\su{\mathfrak{su}}

\def\X{\mathfrak{X}}

\def\D{\mathcal{D}}

\def\C{\mathbb{C}}
\def\R{\mathbb{R}}

\def\N{\mathbb{N}}

				\def\ad{\operatorname{ad}}
				
				\def\Tr{\operatorname{Tr}}
				\def\Id{\operatorname{Id}}
				\def\Ric{\operatorname{Ric}}

				\def\alt{\raise1pt\hbox{$\bigwedge$}}

				\def\del{\partial}
				\def\delbar{\bar{\partial}}
				
				\theoremstyle{plain}
				\newtheorem{theorem}{\bf Theorem}[section]
				\newtheorem{corollary}[theorem]{\bf Corollary}
				\newtheorem{proposition}[theorem]{\bf Proposition}
				\newtheorem{lemma}[theorem]{\bf Lemma}
				
				\theoremstyle{definition}
				
				\newtheorem{example}[theorem]{\bf Example}
				\theoremstyle{remark}
				\newtheorem{remark}[theorem]{Remark}
				
				\newcommand{\ri}{{\rm (i)}}
				\newcommand{\rii}{{\rm (ii)}}
				\newcommand{\riii}{{\rm (iii)}}
				\newcommand{\riv}{{\rm (iv)}}

				\title[Products of Sasakian manifolds]{Harmonic complex structures and special Hermitian metrics on products of Sasakian manifolds}
				
				\author{Adri\'an Andrada}
				\email{adrian.andrada@unc.edu.ar}
				
				\author{Alejandro Tolcachier}
				\email{atolcachier@unc.edu.ar}

				\address{FAMAF, Universidad Nacional de C\'ordoba and CIEM-CONICET, Av. Medina Allende s/n, Ciudad Universitaria, X5000HUA C\'ordoba, Argentina}
				
				\thanks{This work was partially supported by CONICET, SECyT-UNC and FONCyT (Argentina) and the MATHAMSUD Regional Program 21-MATH-06.}
				
				
\begin{document}
					
					\maketitle
					
					\begin{abstract} It is well known that the product of two Sasakian manifolds carries a 2-parameter family of Hermitian structures $(J_{a,b},g_{a,b})$. We show in this article that the complex structure $J_{a,b}$ is harmonic with respect to $g_{a,b}$, i.e. it is a critical point of the Dirichlet energy functional. Furthermore, we also determine when these Hermitian structures are locally conformally Kähler, balanced, strong Kähler with torsion, Gauduchon or $k$-Gauduchon ($k\geq 2$). Finally, we study the Bismut connection associated to $(J_{a,b}, g_{a,b})$ and we provide formulas for the Bismut-Ricci tensor $\Ric^B$ and the Bismut-Ricci form $\rho^B$. We show that these tensors vanish if and only if each Sasakian factor is $\eta$-Einstein with appropriate constants and we also exhibit some examples fulfilling these conditions, thus providing new examples of Calabi-Yau with torsion manifolds. 
					\end{abstract}
					
					\section{Introduction}
					Kähler manifolds are the most celebrated Hermitian manifolds; they lie in the intersection of complex, symplectic, differential and algebraic geometry and due to this interplay they have many interesting geometric properties. A great deal of important complex manifolds are Kähler. Nevertheless, the existence of a Kähler metric imposes strong topological obstructions to the manifold when it is compact, for instance its odd Betti numbers are even. As a consequence, many well-known compact complex manifolds do not admit any Kähler metric; for instance, the Hopf manifolds $\mathbb{S}^1 \times \mathbb{S}^{2n+1}$ with $n\geq 1$. In the search of simply-connected examples, Calabi and Eckmann showed that the manifolds $\mathbb{S}^{2p+1} \times \mathbb{S}^{2q+1}$, with $p,q\geq 1$, carry complex structures but do not admit any Kähler metrics. 
					
					In order to understand non-Kähler Hermitian geometry different approaches have been adopted. One manner is to study conditions on the fundamental 2-form $\omega$ of the Hermitian manifold that are weaker than the Kähler condition $d\omega=0$. For instance,  
					\textit{balanced} metrics are defined by the condition $d\omega^{n-1}=0$, where $n$ is the complex dimension of the manifold, whereas \textit{locally conformally K\"ahler} (LCK) metrics are given by the condition $d\omega=\theta\wedge \omega$, where $\theta$ is a closed $1$-form called the Lee form. Balanced metrics  correspond to the class $\mathcal{W}_3$ in the classification of Gray-Hervella \cite{GH}, and its study started with Michelsohn's work \cite{Mic}. On the other hand, LCK metrics correspond to the class $\mathcal{W}_4$ in the classification of Gray-Hervella \cite{GH}, and they were widely studied since Vaisman's seminal work \cite{Va}. Other conditions which involve the $\del$ and $\delbar$ operators are the following ones: \textit{strong K\"ahler with torsion (SKT)}, given by $\del\delbar\omega=0$; \textit{astheno-Kähler}, given by $\del\delbar \omega^{n-2}=0$; \textit{Gauduchon}, given by $\del\delbar \omega^{n-1}=0$, and its generalization $k$-Gauduchon given by $\del\delbar \omega^{k} \wedge \omega^{n-k-1}=0$.
					
					Some of these non-Kähler Hermitian structures have already been studied on Calabi-Eckmann manifolds. For instance, Michelsohn \cite{Mic} proved that Calabi-Eckmann manifolds $\mathbb{S}^{2p+1}\times \mathbb{S}^{2q+1}$, with $p+q\geq 1$, do not admit any balanced metrics due to homological reasons. On the other hand,
					Wood showed in \cite{Wood} that the Calabi-Eckmann manifolds $\mathbb{S}^{2p+1} \times \mathbb{S}^{2q+1}$ equipped with the product of the round metrics are LCK if and only if $p=0$ or $q=0$. More recently, Cavalcanti established in \cite{Cav} that $\mathbb{S}^1\times \mathbb{S}^3$ and $\mathbb{S}^3 \times \mathbb{S}^3$ are the only Calabi-Eckmann manifolds which admit non-Kähler SKT metrics. We point out that several of these Hermitian properties of the Calabi-Eckmann manifolds were shown to hold using strongly the geometric properties of the spheres (for instance, their (co)homology, the fact that the round metrics have constant sectional curvature,  the existence of the Hopf fibration $\mathbb{S}^{2n+1}\to \C\mathbb{P}^n$, etc.).
					
					Another possible approach is, on a given Riemannian manifold $(M,g)$, to try to detect an optimal almost complex structure from a variational point of view. Assuming that there is at least one  almost complex structure compatible with the metric, the idea is to consider functionals defined on the set of orthogonal almost complex structures on $M$ and then look for extrema of those func\-tio\-nals. A natural functional which has been analyzed in \cite{Wood-Crelle, Wood} is the Dirichlet energy functional $E$, defined as
					\begin{equation}\label{eq:energy}
						E(J):=\int_M \|\nabla J\|^2 \operatorname{vol}_g, 
					\end{equation}
					when $M$ is compact. The first step in the search of (local) minima is to compute the critical points of the functional $E$, which are called \textit{harmonic} almost complex structures. It was shown in \cite{Wood-Crelle} (see also \cite{Wood}) that an orthogonal almost complex structure $J$ is harmonic if and only if 
					\begin{equation*}
						[J,\nabla^*\nabla J]=0, 
					\end{equation*}
					where $\nabla^*\nabla J$ is the \textit{rough Laplacian} of $J$ defined by $\nabla^*\nabla J=\Tr \nabla^2J$. On a non-compact manifold, one can either take the Euler-Lagrange equation $[J,\nabla^*\nabla J]=0$
					as the definition of harmonicity or define the energy on an open subset with compact closure of the manifold and consider variations with compact support included in this subset, and the resulting equation on the open subset is the
					same as in the compact case. 
					
					It follows clearly from \eqref{eq:energy} that K\"ahler structures are harmonic; indeed, they are absolute minimizers. In the more general Hermitian setting, it was shown in \cite{GM} that the complex structure of a balanced Hermitian manifold is harmonic as well as in a LCK manifold provided that the complex dimension is greater than 2. Other examples of harmonic almost complex structures appear on the Calabi-Eckmann manifolds equipped with the product of round metrics \cite{Wood}. Very recently, He and Li introduced in \cite{HL} the \textit{harmonic heat flow} for almost complex structures compatible with a fixed Riemannian metric, which is  a tensor-valued version of the  harmonic map heat equation first studied by Eells-Sampson \cite{ES}. More on harmonic almost complex structures can be found in, for instance, \cite{BLS,DM1,DM2,Lou,LS}. 
					
					A third approach to better comprehend non-K\"ahler Hermitian manifolds is by considering connections that preserve the Hermitian structure, i.e. both the Riemannian metric and the complex structure are parallel. Such a connection is called Hermitian. It is well known that the Levi-Civita connection associated to the Hermitian metric is a Hermitian connection if and only if the manifold is K\"ahler, so in the non-K\"ahler setting we should look for other such  connections. In fact, on any given Hermitian manifold there are infinitely many Hermitian connections, but there is only one satisfying an extra condition: its torsion, considered as a $(0,3)$-tensor by contracting with the metric, is actually a $3$-form \cite{Bi}. This connection, that we denote $\nabla^B$, is called the \textit{Bismut} connection, although in the physics literature it is known as the \textit{K\"ahler with torsion} (or KT) connection, and more recently also the name \textit{Strominger} connection has been used. Since both the complex structure and the Hermitian metric are $\nabla^B$-parallel we have that its holonomy group  $\operatorname{Hol}^B$ is contained in $\operatorname{U}(n)$, where $2n$ is the real dimension of the manifold. In particular, $2n$-dimensional Hermitian manifolds whose (restricted) Bismut holonomy group is contained in $\operatorname{SU}(n)$ have attracted plenty of attention (see for instance \cite{AV,FG,GGP,Grant,IP,UV}). These manifolds are known as \textit{Calabi-Yau with torsion} (CYT for short), and they appear in heterotic string theory, related to the Strominger system in six dimensions (see for instance \cite{H,St}). The CYT condition is equivalent to the vanishing of the Bismut-Ricci form (see \eqref{eq:Ric-rho} below for its definition). We point out that Bismut flat manifolds (i.e., the Bismut curvature tensor $R^B$ vanishes) have recently been characterized in \cite{WYZ}: if $M$ is a compact Hermitian manifold with flat Bismut connection, then its universal cover is a Lie group $G$ equipped with a bi-invariant metric and a left invariant complex structure compatible with the metric. In particular, $G$ is the product of a compact semisimple Lie group and a real vector space. 
					
					Our first goal in this article is to generalize some of the properties of Calabi-Eckmann manifolds mentioned above to the product of two arbitrary Sasakian manifolds, since it is well known that odd-dimensional spheres carry a canonical Sasakian structure. It was shown by Morimoto \cite{Mor} that the product of two normal almost contact manifolds has a natural complex structure. This was later generalized independently by Tsukada \cite{Tsu} and Watson \cite{Wat}, who showed the existence of a family of complex structures $J_{a,b}$ for $a,b\in\R$, $b\neq 0$, which correspond to the complex structures on $\mathbb{S}^{2p+1}\times \mathbb{S}^{2q+1}$ given in \cite{CE}. They also showed the existence of a family of compatible Hermitian metrics $g_{a,b}$. We will restrict to the case of a product of two Sasakian manifolds and the corresponding Hermitian structures $(J_{a,b}, g_{a,b})$ will be the central object of study throughout the paper. 
					
					A second goal of this article is to study the Bismut connection associated to the Hermitian structure $(J_{a,b},g_{a,b})$. Concretely, we study the vanishing of the Bismut-Ricci tensor $\Ric^B$ and the Bismut-Ricci form $\rho^B$. It will turn out that these conditions are closely related to a particular family of Sasakian manifolds, called \textit{$\eta$-Einstein}. A Sasakian manifold is called $\eta$-Einstein if the Ricci tensor of the Sasakian metric satisfies $\Ric=\lambda g+\nu \eta\otimes \eta$ for certain constants $\lambda, \nu\in\R$, where $\eta$ is the $1$-form dual to the Reeb vector field. 
					
					The article is structured as follows. In \S\ref{sec:prelim} we recall basic notions on Sasakian manifolds and their transverse geometry and present some preliminary results. In \S\ref{sec:prod-sasakian} we study the Levi-Civita connection of the metric $g_{a,b}$ on the product $S_1\times S_2$, where $S_1$ and $S_2$ are Sasakian manifolds, and we use this in \S\ref{sec:harmonic} to show that $J_{a,b}$ is harmonic on $(S_1\times S_2, g_{a,b})$  (see Theorem \ref{theorem:prod-sasakian}). Next, in \S\ref{sec:Hermitian}, we study the balanced, LCK, SKT and $k$-Gauduchon ($k\geq 2$) conditions on $S_1 \times S_2$. We show in Theorem \ref{theorem:k-gau} that the Hermitian structure $(J_{a,b},g_{a,b})$ on $S_1\times S_2$ is always Gauduchon (i.e. $(n-1)$-Gauduchon) and it is $k$-Gauduchon ($2\leq k\leq n-2)$ if and only if is astheno-Kähler. This complements the result in \cite{FU}, where it was shown that $(J_{a,b},g_{a,b})$ is 1-Gauduchon if and only if it is astheno-Kähler. The astheno-Kähler condition was previously characterized in \cite{Mat}. Finally, in \S\ref{sec:Bismut} we provide an explicit expression for the Bismut connection associated to $(J_{a,b}, g_{a,b})$ in terms of the characteristic connections on $S_1$ and $S_2$. As an application of this explicit expression, we provide formulas for the Bismut-Ricci tensor $\Ric^B$ and the Bismut-Ricci form $\rho^B$, and determine when they vanish (Theorems \ref{theorem:BRF} and \ref{theorem:CYT}, respectively). More precisely, we show that $\Ric^B=0$ or $\rho^B=0$ hold if and only if both Sasakian factors $S_1$ and $S_2$ are $\eta$-Einstein with certain appropriate constants $(\lambda_1,\nu_1)$ and $(\lambda_2,\nu_2)$; and we exhibit examples of Hermitian structures $(J_{a,b},g_{a,b})$ with $\Ric^B=0$ or $\rho^B=0$, leading to new examples of CYT structures.
					
					\medskip 
					
					\textbf{Acknowledgments.} 
					The authors are grateful to Romina Arroyo, Jorge Lauret, Henrique Sá Earp, Mauro Subils and Jeffrey Streets for their useful comments and suggestions. The authors would also like to thank the hospitality of the Instituto de Matemática, Estatística e Computação Científica at UNICAMP (Brazil), where they were introduced to the theory of harmonic almost complex structures. 
					
					\medskip
					
					\section{Preliminaries on Sasakian manifolds}\label{sec:prelim}
					
					An \textit{almost contact structure} on a differentiable manifold $M^{2n+1}$ is a triple $(\varphi,\xi,\eta)$, where $\varphi$ is a (1,1)-type tensor field, $\xi$ a vector field, and $\eta$ a $1$-form satisfying
					\begin{equation*}
						\varphi^2=-\Id+\eta\otimes \xi,\quad \eta(\xi)=1,
					\end{equation*}  
					which imply $\varphi(\xi)=0$ and $\eta\circ\varphi=0$. $(M,\varphi,\xi,\eta)$ is called an \textit{almost contact manifold} and $\xi$ is called the Reeb vector field. The tangent bundle of $M$ splits as $TM={\mathcal D}\oplus{\mathcal L}$, where $\mathcal D=\mathrm{Ker}\,\eta=\mathrm{Im}\,\varphi$ and ${\mathcal L}$ is the line bundle spanned by $\xi$.
					
					On the product manifold $M\times \mathbb{R}$ there is a natural almost complex structure $J$ defined by
					\begin{equation}\label{eq:J-en-MxR}
						J\left(X+f\frac{d}{dt}\right)=\varphi
						X-f\xi+\eta(X)\frac{d}{dt},
					\end{equation}
					where $X\in\X(M)$, $t$ is the coordinate on $\mathbb{R}$ and $f$ is a smooth function on $M\times \mathbb{R}$. If $J$ is integrable, the almost contact structure is said to be \emph{normal}. This is equivalent to the
					vanishing of the tensor field
					\begin{equation*}
						N_\varphi:=[\varphi,\varphi]+d\eta\otimes\xi,
					\end{equation*}
					where $[\varphi,\varphi]$ is the Nijenhuis torsion of $\varphi$ defined by
					\[ [\varphi,\varphi](X,Y)=[\varphi X,\varphi Y]+\varphi^2[X,Y]-\varphi[\varphi X,Y]-\varphi[X,\varphi Y].\]
					Note that if the almost contact structure is normal then
					\begin{equation}\label{eq:fi-xi}
						\varphi[\xi,X]=[\xi,\varphi X] \qquad \text{for all } X\in\mathfrak{X}(M).
					\end{equation}
					In particular, 
					\begin{equation}\label{eq:eta-ad_xi}
						[\xi,X]\in \Gamma(\D) \qquad \text{for all } X\in \Gamma(\mathcal D).
					\end{equation}
					
					\medskip
					
					An \textit{almost contact metric structure} on $M$ is $(\varphi,\xi,\eta, g)$, where $(\varphi,\xi,\eta)$ is an almost contact structure and $g$ is a Riemannian metric on $M$ satisfying
					\begin{equation}\label{eq:metric}
						g(\varphi X,\varphi Y)=g(X,Y)-\eta(X)\eta(Y), \quad X,Y\in \mathfrak{X}(M).
					\end{equation} 
					This equation implies:
					\[  g(\varphi X, Y)=-g(X,\varphi Y)\qquad \text{and} \qquad g(\xi,X)=\eta(X),\] for all $X,Y\in\X(M)$.
					That is, $\varphi$ is skew-symmetric and the vector field $\xi$ is $g$-dual to the $1$-form $\eta$. In analogy with the almost Hermitian setting, the fundamental $2$-form $\Phi$ can be defined by
					\[ \Phi (X,Y)=g(X,\varphi Y).\]
					
					A normal almost contact manifold $S$ is called \textit{Sasakian} if $d\eta= 2 \Phi$. In particular, $\Phi$ is exact and $\eta$ is a contact form on $S$, i.e. $\eta\wedge (d\eta)^n$ is a volume form. Sasakian manifolds form the most important class of almost contact metric manifolds, due to their close relation with K\"ahler manifolds. Indeed, the Riemannian cone of an almost contact metric manifold endowed with the almost complex structure given by \eqref{eq:J-en-MxR} is K\"ahler if and only if the structure is Sasakian.
					
					\smallskip
					
					Some properties of Sasakian manifolds that we will need in forthcoming sections are stated in the following lemma.
					
					\begin{lemma}\label{lemma:propiedades}
						If $(S,\varphi,\xi,\eta, g)$ is a Sasakian manifold with fundamental $2$-form $\Phi(X,Y)=g(X,\varphi Y)$ and $d\eta=2\Phi$, then:
						\begin{enumerate}
							\item[$\ri$] $\xi$ is a unit Killing vector field on $S$,
							\item[$\rii$] $\nabla_X \xi=-\varphi X$ for all $X\in \mathfrak{X}(S)$; in particular, $\nabla_\xi\xi=0$,
							\item[$\riii$] $\nabla_\xi X=[\xi,X]-\varphi X$ for all $X\in \mathfrak{X}(S)$,
							\item[$\riv$] $(\nabla_X \varphi)Y=g(X,Y)\xi-\eta(Y) X$ for all $X,Y\in \mathfrak{X}(S)$.
						\end{enumerate}
					\end{lemma}
					
					The proof of Lemma \ref{lemma:propiedades} is standard and can be found in \cite{Blair}. 
					
						
					
					\medskip
					
					\subsection{Transverse geometry of Sasakian manifolds}
					
					Let $(S,\varphi,\eta,\xi,g)$ be a Sasakian manifold of dimension $2n+1$. Recall that $\D=\operatorname{Ker} \eta=\operatorname{Im}\varphi$ is a subbundle of $TS$ of rank $2n$. The Sasakian structure induces on $\D$ a natural connection $\nabla^T$, called the transverse Levi-Civita connection which, for any $U\in \Gamma(\D)$, is defined by
					\begin{equation}\label{eq:transverse-0}
						\nabla^T_\xi U= [\xi,U], \qquad \nabla^T_X U= (\nabla_X U)^\D, \quad X\in  \Gamma(\D),
					\end{equation}
					where $(\cdot)^\D$ denotes the projection onto $\D$. This is the only connection on $\mathcal D$ that satisfies
					\begin{equation}\label{eq:transverse-1}
						\nabla^T_X(\varphi|_\D)=0, \qquad \nabla^T_X (g|_\D)=0, \qquad \nabla^T_U V-\nabla^T_V U=[U,V]^\D,
					\end{equation} 
					for any $X\in \X(S)$ and $U,V\in\Gamma(\D)$. Note that 
					\begin{equation}\label{eq:transverse-2}
						\nabla_UV=-\Phi(U,V)\xi+\nabla_U^TV, \quad U,V\in \Gamma(\D),
					\end{equation}
					which implies
					\begin{equation}\label{eq:transverse-3}
						[U,V]=-2\Phi(U,V)\xi+[U,V]^\mathcal{D}, \quad U,V\in \Gamma(\D).
					\end{equation} 
					Using \eqref{eq:transverse-0}--\eqref{eq:transverse-3} we obtain the following result.
					
					\begin{lemma}\label{lemma:transverse}
						For any $U,V,W\in\Gamma(\D)$ the following identities hold:
						\begin{enumerate}
							\item[$\ri$] $\nabla^T_{[U,V]^\mathcal{D}}W=\nabla^T_{[U,V]}W+2\Phi(U,V)[\xi,W]$,
							\item[$\rii$]
							$\nabla_{[U,V]} W=2 \Phi(U,V)\varphi W-\Phi([U,V]^{\mathcal{D}},W)\xi+\nabla^T_{[U,V]} W$,
							\item[$\riii$] $R(U,V)W=R^T(U,V)W+\Phi(V,W)\varphi U -\Phi(U,W)\varphi V-2\Phi(U,V)\varphi W$, 
							\item[$\riv$] $R(U,V)\xi =0$.
						\end{enumerate}
					\end{lemma}
					
					\begin{proof}
						(i) is a straightforward consequence of \eqref{eq:transverse-3} and \eqref{eq:transverse-0}:
						\begin{align*}
							\nabla^T_{[U,V]}W & = \nabla^T_{(-2\Phi(U,V)\xi+[U,V]^\mathcal{D})}W \\
							& = -2\Phi(U,V)[\xi,W]+\nabla^T_{[U,V]^\mathcal{D}}W.
						\end{align*}
						
						For (ii), using \eqref{eq:transverse-3} and \eqref{eq:transverse-2} we compute
						\begin{align*}
							\nabla_{[U,V]} W&= \nabla_{(-2\Phi(U,V)\xi+[U,V]^{\mathcal{D}})} W\\
							&= -2\Phi(U,V) ([\xi,W]-\varphi W)+\nabla_{[U,V]^{\mathcal{D}}} W\\
							&= -2\Phi(U,V)([\xi,W]-\varphi W)-\Phi([U,V]^{\mathcal{D}},W)\xi+\nabla^T_{[U,V]^{\mathcal{D}}} W\\
							&=-2\Phi(U,V)([\xi,W]-\varphi W)-\Phi([U,V]^{\mathcal{D}},W)\xi+\nabla^T_{[U,V]} W\\
							&\quad +2\Phi(U,V)[\xi,W]\\
							&=2\Phi(U,V)\varphi W-\Phi([U,V]^{\mathcal{D}},W)\xi+\nabla^T_{[U,V]}W,
						\end{align*}
						where we have used (i) in the fourth equality.
						
						For (iii), beginning with the definition $R(U,V)W=\nabla_U\nabla_VW-\nabla_V\nabla_U W-\nabla_{[U,V]}W$ and using \eqref{eq:transverse-0}, \eqref{eq:transverse-3}, Lemma \ref{lemma:propiedades} and (ii) we obtain
						\begin{align*}
							R(U,V)W & = -U(\Phi(V,W))\xi +\Phi(V,W)\varphi U-\Phi(U,\nabla^T_VW)\xi +\nabla^T_U\nabla^T_V W\\
							& \quad +V(\Phi(U,W))\xi-\Phi(U,W)\varphi V+\Phi(V,\nabla^T_UW)\xi -\nabla^T_V\nabla^T_U W \\
							& \quad  -(2\Phi(U,V)\varphi W-\Phi([U,V]^{\mathcal{D}},W)\xi+\nabla^T_{[U,V]} W) 
						\end{align*}
						Next, using (i), the definition of $\Phi$ and \eqref{eq:transverse-1} we arrive at
						\begin{align*}
							R(U,V)W & = -g(\nabla^T_UV,\varphi W)\xi - g(V,\varphi \nabla^T_U W)\xi +\Phi(V,W)\varphi U -g(U, \varphi \nabla^T_V W)\xi \\ 
							& \quad + g(\nabla^T_VU,\varphi W)\xi+g(U,\varphi \nabla^T_VW)\xi -\Phi(U,W)\varphi V + g(V,\varphi \nabla^T_UW)\xi\\
							& \quad -2\Phi(U,V)\varphi W+\Phi([U,V]^{\mathcal{D}},W)\xi+ R^T(U,V)W\\
							& = R^T(U,V)W+\Phi(V,W)\varphi U -\Phi(U,W)\varphi V-2\Phi(U,V)\varphi W,
						\end{align*}
						and (iii) is proved.
						
						For (iv), we compute 
						\begin{align*}
							R(U,V)\xi & = \nabla_U \nabla_{V} \xi-\nabla_{V} \nabla_U \xi-\nabla_{[U,V]}\xi \\
							& = -\nabla_U \varphi V+\nabla_{V}\varphi U+\varphi[U,V]\\
							& = \Phi(U,\varphi V)\xi-\nabla^T_U \varphi V-\Phi(V,\varphi U) \xi+\nabla_V^T \varphi U+\varphi [U,V]^\mathcal{D}\\
							& = -g(U,V)\xi-\varphi \nabla_U^T V+g(V,U)\xi+\varphi \nabla^T_V U+\varphi[U,V]^{\mathcal{D}},
						\end{align*}
						according to Lemma \ref{lemma:propiedades}$\rii$. It follows from \eqref{eq:transverse-1} that this last expression vanishes, therefore $R(U,V)\xi=0$, and the proof is complete.
					\end{proof}
					
					\medskip
					
					\begin{corollary}\label{corollary:R-commutes}
						For any $U,V\in \Gamma(\D)$, each curvature endomorphism $R(U,V)$ preserves $\D$. Moreover, for $V=\varphi U$ we have that $R(U,\varphi U)|_\D$ commutes with $\varphi|_\D$.
					\end{corollary}
					
					\medskip 
					
					\section{Hermitian structures on the product of Sasakian manifolds} \label{sec:prod-sasakian}
					
					We recall next the following construction, developed independently by Tsukada \cite{Tsu} and Watson \cite{Wat}, both based on a previous construction due to Morimoto \cite{Mor}, using ideas from \cite{CE}. With this construction, one can define a Hermitian structure on the product of two manifolds equipped with normal almost contact metric structures. We will focus later on the product of Sasakian manifolds. 
					
					Let $M_1$ and $M_2$ be differentiable manifolds of dimension $2n_1+1$ and $2n_2+1$ and let $(\varphi_1,\xi_1,\eta_1,g_1)$ and $(\varphi_2,\xi_2,\eta_2,g_2)$ be almost contact metric structures on $M_1$ and $M_2$, respectively. 
					
					For $a,b\in\R, \, b\neq 0$, we can induce an almost Hermitian structure $(J_{a,b},g_{a,b})$ on the product manifold $M:=M_1\times M_2$ as follows: for $X_1\in \mathfrak{X}(M_1)$ and $X_2\in \mathfrak{X}(M_2)$, define an almost complex structure $J_{a,b}$ on $M$ by
					\begin{align} \label{eq:Jab}
						J_{a,b} (X_1+X_2) & =  \varphi_1 X_1- \left(\frac{a}{b}\eta_1(X_1)+\frac{a^2+b^2}{b}\eta_2(X_2)\right)\xi_1   \\  
						& \quad +\varphi_2 X_2+\left(\frac{1}{b}\eta_1(X_1)+\frac{a}{b}\eta_2(X_2)\right)\xi_2.\nonumber
					\end{align}
					Next, define a Riemannian metric $g_{a,b}$ on $M$ by 
					\begin{align}\label{eq:gab}
						g_{a,b} (X_1+X_2,Y_1+Y_2) & = g_1(X_1,Y_1)+a[\eta_1(X_1)\eta_2(Y_2)+\eta_1(Y_1)\eta_2(X_2)] \\
						& \quad +g_2(X_2,Y_2)+(a^2+b^2-1)\eta_2(X_2)\eta_2(Y_2). \nonumber
					\end{align}
					It is an easy exercise on quadratic forms to verify that $g_{a,b}$ is indeed positive definite and $J_{a,b}$ is Hermitian with respect to $g_{a,b}$.
					
					Regarding $\mathfrak{X}(M_1)$ and $\mathfrak{X}(M_2)$ as subalgebras of $\mathfrak{X}(M)$ in a natural manner, \eqref{eq:Jab} and \eqref{eq:gab} can be rewritten in the following way, where $U_i\in\Gamma(\mathcal D_i)$:
					\begin{gather*}
						J_{a,b}\xi_1   =  -\frac{a}{b}\xi_1+\frac{1}{b}\xi_2, \qquad   J_{a,b}U_1  = \varphi_1 U_1, \\
						J_{a,b}\xi_2  = -\frac{a^2+b^2}{b}\xi_1+\frac{a}{b}\xi_2, \qquad  J_{a,b}U_2  = \varphi_2 U_2,  
					\end{gather*}
					and, for $X_i,Y_i\in \mathfrak{X}(M_i)$:
					\begin{align*} 
						g_{a,b}(X_1,Y_1) &= g_1(X_1,Y_1),\qquad   g_{a,b}(X_1,X_2)= a\eta_1(X_1)\eta_2(X_2)  \\
						g_{a,b}(X_2,Y_2) &= g_2(X_2,Y_2)+(a^2+b^2-1)\eta_2(X_2)\eta_2(Y_2), 
					\end{align*}
					Note that $g_{a,b}$ coincides with $g_1$ on $M_1$ and with $g_2$ on $\D_2$, but it  modifies the length of $\xi_2$ by a factor of $a^2+b^2$; also, $\xi_1$ and $\xi_2$ are no longer orthogonal whenever $a\neq 0$. Moreover, $g_{a,b}$ is the product $g_1\times g_2$ if and only if $a=0, b=\pm 1$.
					
					Morimoto's original construction corresponds to the case $a=0$, $b=1$. He proved the following result:
					
					\begin{proposition}\cite[Proposition 3]{Mor}
						Let $(\varphi_1,\xi_1,\eta_1)$ and $(\varphi_2,\xi_2,\eta_2)$ be almost contact structures on $M_1$ and $M_2$, respectively. Then the almost complex structure $J_{0,1}$ on $M=M_1\times M_2$ is integrable if and only if both almost contact structures are normal.
					\end{proposition}
					
					More generally, the following result can be proved in the same way as in \cite{Mor}:
					
					\begin{proposition}
						Let $(\varphi_1,\xi_1,\eta_1)$ and $(\varphi_2,\xi_2,\eta_2)$ be almost contact structures on $M_1$ and $M_2$, respectively. If both almost contact structures are normal then the almost complex structure $J_{a,b}$ is integrable for any $a\in\R,\, b\in\R,\, b\neq 0$.
					\end{proposition}
					
					From now on, we will deal only with the case when $(S_1,\varphi_1,\xi_1,\eta_1,g_1)$ and $(S_2,\varphi_2,\xi_2,\eta_2,g_2)$ are Sasakian manifolds. We will denote $M_{a,b}:=S_1\times S_2$ equipped with the Hermitian structure $(J_{a,b},g_{a,b})$. Moreover, we will denote simply $J:=J_{a,b}$, $g:=g_{a,b}$ since there will be no risk of confusion. 
					
					In forthcoming sections we will need explicit formulas for the Levi-Civita connection $\nabla$ on $M_{a,b}$ associated to $g$ in terms of the Levi-Civita connections $\nabla^1$ and $\nabla^2$ on $(S_1,g_1)$ and $(S_2,g_2)$, respectively. We will use the following expressions which appear for instance in \cite{LPS}: given $X_i,Y_i,Z_i\in \X(S_i)$, we have that
					\begin{align}
						g(\nabla_{X_1} Y_1, Z_1)&=g_1(\nabla^1_{X_1} Y_1, Z_1), \quad g(\nabla_{X_1} Y_1, Z_2)=a \eta_1 (\nabla^1_{X_1} Y_1) \eta_2(Z_2) \nonumber\\ 
						g(\nabla_{X_2} Y_2, Z_1)&=a \eta_2(\nabla^2_{X_2} Y_2) \eta_1(Z_1) \nonumber\\  
						g(\nabla_{X_2} Y_2, Z_2)&=g_2(\nabla^2_{X_2} Y_2,Z_2)+(a^2+b^2-1)[\eta_2(\nabla^2_{X_2}Y_2)\eta_2(Z_2) \label{eq:nabla}\\
						&\quad -\eta_2(X_2)g_2(\varphi_2 Y_2,Z_2)-\eta_2(Y_2)g_2(\varphi_2 X_2,Z_2)] \nonumber \\
						g(\nabla_{X_1} Y_2, Z_1)&=-a \eta_2(Y_2)g_1(\varphi_1 X_1,Z_1), \quad 
						g(\nabla_{X_1} Y_2, Z_2)=-a \eta_1(X_1)g_2(\varphi_2 Y_2,Z_2) \nonumber\\
						g(\nabla_{X_2} Y_1, Z_1)&=-a \eta_2(X_2)g_1(\varphi_1 Y_1,Z_1), \quad 
						g(\nabla_{X_2} Y_1, Z_2)=-a \eta_1(Y_1) g_2(\varphi_2 X_2,Z_2) \nonumber
					\end{align}
					
					\medskip
					
					The next result follows from the set of equations \eqref{eq:nabla}:
					
					\begin{corollary}\label{corollary:nabla} With notation as above,
						\begin{enumerate}
							\item[$\ri$] $\nabla_{X_1} Y_1=\nabla^1_{X_1} Y_1\in \X(S_1)$,
							\item[$\rii$] $\nabla_{X_2}Y_2=\nabla^2_{X_2}Y_2-(a^2+b^2-1)[\eta_2(X_2)\varphi_2Y_2+\eta_2(Y_2)\varphi_2 X_2]\in\X(S_2)$,
							\item[$\riii$] $\nabla_{X_1}Y_2=-a[\eta_2(Y_2)\varphi_1X_1+\eta_1(X_1)\varphi_2Y_2]\in \X(S_1)\oplus \X(S_2)$,
							\item[$\riv$] $\nabla_{X_2}Y_1=-a[\eta_2(X_2)\varphi_1Y_1+\eta_1(Y_1)\varphi_2X_2]\in \X(S_1)\oplus \X(S_2)$.
						\end{enumerate}
						In particular, $\nabla_{\xi_1}\xi_1=\nabla_{\xi_2}\xi_2=\nabla_{\xi_1} \xi_2=\nabla_{\xi_2} \xi_1=0$.
					\end{corollary}
					
					\medskip
					
					Using the previous corollary we compute next $\nabla J$, which will be needed in the proof of Lemma \ref{lemma:codif-J} below.
					
					\begin{lemma}\label{lemma:nabla-XX}
						For any $X_i,Y_i\in \mathfrak{X}(S_i)$, $i=1,2$, 
						\begin{enumerate}
							\item[$\ri$] $(\nabla_{X_1}J)Y_1=g_1(X_1,Y_1)\xi_1-\eta_1(Y_1)X_1-\frac{a}{b}\Phi_1(X_1,Y_1)\xi_1+\frac{1}{b}\Phi_1(X_1,Y_1)\xi_2$,
							\item[$\rii$] 
							$\begin{aligned}[t] (\nabla_{X_2}J)Y_2&=[g_2(X_2,Y_2)+(a^2+b^2-1)\eta_2(X_2)\eta_2(Y_2)]\xi_2-(a^2+b^2)\eta_2(Y_2)X_2\\
								&\quad-\frac{a^2+b^2}{b}\Phi_2(X_2,Y_2)\xi_1+\frac{a}{b}\Phi_2(X_2,Y_2)\xi_2
							\end{aligned}$
							\item[$\riii$] $(\nabla_{X_1} J)Y_2=a\eta_2(Y_2)\eta_1(X_1)\xi_1-a\eta_2(Y_2)X_1+b\eta_2(Y_2)\varphi_1 X_1$.
							\item[$\riv$] $(\nabla_{X_2} J)Y_1= a[\eta_1(Y_1)\eta_2(X_2)\xi_2-\eta_1(Y_1)X_2]-b\eta_1(Y_1)\varphi_2 X_2$.
						\end{enumerate}
						In particular, $\nabla_{\xi_1}J=0$ and $\nabla_{\xi_2}J=0$.
					\end{lemma}
					
					\begin{proof}
						We compute $\nabla J$ using Corollary \ref{corollary:nabla} and the definition of $J$.
						
						For (i), 
						\begin{align*}
							(\nabla_{X_1} J)Y_1 &=\nabla_{X_1} JY_1-J\nabla^1_{X_1} Y_1
						\end{align*} 
						We will expand each term in detail. On the one hand, 
						\begin{align*}  \nabla_{X_1} JY_1&=\nabla_{X_1}(\varphi_1 Y_1-\frac{a}{b} \eta_1(Y_1)\xi_1+\frac{1}{b}\eta_1(Y_1)\xi_2)\\
							&=\nabla_{X_1}^1 \varphi_1 Y_1-\frac{a}{b} (X_1(\eta_1(Y_1)) \xi_1-\eta_1(Y_1)\varphi_1 X_1)\\
							&\quad +\frac{1}{b}(X_1(\eta_1(Y_1))\xi_2-a \eta_1(Y_1)\varphi_1 X_1).
						\end{align*} On the other hand,
						\[ -J\nabla^1_{X_1} Y_1=-\varphi_1 \nabla^1_{X_1} Y_1+\frac{a}{b} \eta_1(\nabla^1_{X_1} Y_1)\xi_1-\frac{1}{b} \eta_1(\nabla^1_{X_1} Y_1)\xi_2.\]
						Putting these two expressions together we arrive at
						\begin{align*}
							(\nabla_{X_1} J) Y_1&
							= (\nabla^1_{X_1} \varphi_1) Y_1-\frac{a}{b} g_1(\nabla_{X_1}^1 \xi_1, Y_1)\xi_1+\frac{1}{b} g_1(\nabla^1_{X_1} \xi_1, Y_1)\xi_2\\
							&=g_1(X_1,Y_1)\xi_1-\eta_1(Y_1)X_1-\frac{a}{b} \Phi_1(X_1,Y_1)\xi_1+\frac{1}{b} \Phi_1(X_1,Y_1)\xi_2,
						\end{align*}
						where we have used Lemma \ref{lemma:propiedades}(iv).
						
						Next, for (ii), 
						\begin{align*}
							(\nabla_{X_2} J)Y_2&= \nabla_{X_2} JY_2-J(\nabla^2_{X_2}Y_2-(a^2+b^2-1)[\eta_2(X_2)\varphi_2Y_2+\eta_2(Y_2)\varphi_2 X_2]).
						\end{align*}
						The first term is equal to
						\begin{align*} 
							\nabla_{X_2} JY_2&=\nabla_{X_2}(\varphi_2 Y_2-\frac{a^2+b^2}{b} \eta_2(Y_2)\xi_1+\frac{a}{b} \eta_2(Y_2)\xi_2)\\
							&=\nabla^2_{X_2} \varphi_2 Y_2-(a^2+b^2-1)\eta_2(X_2)\varphi_2^2Y_2\\
							&\quad -\frac{a^2+b^2}{b} (X_2(\eta_2(Y_2))\xi_1-a \eta_2(Y_2)\varphi_2 X_2)\\
							&\quad +\frac{a}{b} (X_2(\eta_2(Y_2))\xi_2-(a^2+b^2)\eta_2(X_2)\varphi_2 X_2),
						\end{align*} 
						and the second term is equal to
						\begin{align*}
							&-\varphi_2(\nabla^2_{X_2}Y_2-(a^2+b^2-1)[\eta_2(X_2)\varphi_2Y_2+\eta_2(Y_2)\varphi_2 X_2])\\
							&\quad +\frac{a^2+b^2}{b}\eta_2(\nabla^2_{X_2}Y_2)\xi_1-\frac{a}{b}\eta_2(\nabla^2_{X_2}Y_2)\xi_2. 
						\end{align*}
						Therefore,
						\begin{align*}
							(\nabla_{X_2} J)Y_2&=(\nabla^2_{X_2} \varphi_2)Y_2+(a^2+b^2-1)\eta_2(Y_2)\varphi_2^2 X_2\\
							&\quad -\frac{a^2+b^2}{b}\Phi_2(X_2,Y_2)\xi_1+\frac{a}{b}\Phi_2(X_2,Y_2)\xi_2\\
							&=(g_2(X_2,Y_2)+(a^2+b^2-1)\eta_2(X_2)\eta_2(Y_2))\xi_2-(a^2+b^2)\eta_2(Y_2)X_2\\
							&\quad-\frac{a^2+b^2}{b}\Phi_2(X_2,Y_2)\xi_1+\frac{a}{b}\Phi_2(X_2,Y_2)\xi_2,
						\end{align*}
						using again Lemma \ref{lemma:propiedades}(iv). 
						
						Now, for (iii) and (iv) we compute
						\begin{align*}
							(\nabla_{X_1} J)Y_2&=\nabla_{X_1}(\varphi_2 Y_2-\frac{a^2+b^2}{b}\eta_2(Y_2)\xi_1+\frac{a}{b}\eta_2(Y_2)\xi_2)\\
							&\quad +a J(\eta_2(Y_2)\varphi_1 X_1+\eta_1(X_1)\varphi_2 Y_2)\\
							&=-a \eta_1(X_1)\varphi_2^2 Y_2+b\eta_2(Y_2)\varphi_1 X_1+a (\eta_2(Y_2)\varphi_1^2 X_1+\eta_1(X_1)\varphi_2^2 Y_2)\\
							&=a\eta_2(Y_2) \varphi_1^2 X_1+b \eta_2(Y_2) \varphi_1 X_1\\
							&=a\eta_2(Y_2)\eta_1(X_1)\xi_1-a\eta_2(Y_2)X_1+b\eta_2(Y_2)\varphi_1 X_1,
						\end{align*}
						\begin{align*}
							(\nabla_{X_2} J)Y_1&=\nabla_{X_2} (\varphi_1 Y_1-\frac{a}{b}\eta_1(Y_1)\xi_1+\frac{1}{b}\eta_1(Y_1)\xi_2)+aJ(\eta_2(X_2)\varphi_1 Y_1+\eta_1(Y_1)\varphi_2 X_2)\\
							&=-a \eta_2(X_2)\varphi_1^2 Y_1-\frac{a}{b} \eta_1(Y_1) (-a\varphi_2 X_2) -\frac{a^2+b^2}{b} \eta_1(Y_1) \varphi_2 X_2\\
							&\quad +a(\eta_2(X_2)\varphi_1^2 Y_1+\eta_1(Y_1)\varphi_2^2 X_2)\\
							&= a[\eta_1(Y_1)\eta_2(X_2)\xi_2-\eta_1(Y_1)X_2]-b\eta_1(Y_1)\varphi_2 X_2.
						\end{align*}
						The last statement follows easily from the previous computations.
					\end{proof}
					
					In the next result $R$ denotes the curvature tensor of $\nabla$, while $R^i$ denotes the curvature tensor of $\nabla^i$, $i=1,2$. We shall compute only the curvature tensors which we will need in the proof of Theorem \ref{theorem:prod-sasakian}. Let us set $\lambda_{a,b}:=a^2+b^2-1$ to shorten a little bit the statement and proof of the lemma.
					
					\begin{lemma}\label{lemma:curvatura} With notation as above, for $U_i,V_i\in \Gamma(\D_i), Z_i\in \X(S_i)$,
						\begin{enumerate}
							\item[$\ri$] $R(\xi_1,\xi_2)=0$,
							\item[$\rii$] $R(U_1,V_1)Z_1=R^1(U_1,V_1)Z_1$ and $R(U_1,V_1)Z_2=-2a \Phi_1(U_1,V_1)\varphi_2 Z_2$,
							\item[$\riii$] $R(U_2,V_2)Z_1=-2 a\Phi_2(U_2,V_2) \varphi_1 Z_1$ and \\ \quad 
							$\begin{aligned}[t] 
								R(U_2,V_2)Z_2&=R^2(U_2,V_2)Z_2+\lambda_{a,b}[\Phi_2(V_2,Z_2)\varphi_2 U_2-\Phi_2(U_2,Z_2)\varphi_2 V_2\\&\quad -2\Phi_2(U_2,V_2)\varphi_2 Z_2].
							\end{aligned}$
						\end{enumerate}
						In particular, $R(U_i,V_i)\xi_1=R(U_i,V_i)\xi_2=0$.
					\end{lemma}
					
					\begin{proof}
						For (i), we compute $R(\xi_1,\xi_2)Z_i$ for $Z_i\in \X(S_i)$, $i=1,2$, using Corollary \ref{corollary:nabla} and properties of Sasakian manifolds such as Lemma \ref{lemma:propiedades} and \eqref{eq:fi-xi}.
						
						For $Z_1$ we have
						\begin{align*}
							R(\xi_1,\xi_2)Z_1 & = \nabla_{\xi_1}\nabla_{\xi_2}Z_1-\nabla_{\xi_2}\nabla_{\xi_1}Z_1 \\
							& = \nabla_{\xi_1}(-a\varphi_1 Z_1)-\nabla_{\xi_2}([\xi_1,Z_1]-\varphi_1 Z_1)\\
							& = -a([\xi_1,\varphi_1 Z_1]-\varphi_1^2 Z_1)+a\varphi_1([\xi_1,Z_1]-\varphi_1 Z_1)\\
							& =0,
						\end{align*}
						and for $Z_2$ we have
						\begin{align*}
							R(\xi_1,\xi_2)Z_2 & = \nabla_{\xi_1}\nabla_{\xi_2}Z_2-\nabla_{\xi_2}\nabla_{\xi_1}Z_2 \\
							& = \nabla_{\xi_1}([\xi_2,Z_2]-\varphi_2 Z_2-(a^2+b^2-1)\varphi_2 Z_2)-\nabla_{\xi_2}(-a\varphi_2 Z_2)\\
							& = -a(\varphi_2[\xi_2,Z_2] -(a^2+b^2) \varphi_2^2 Z_2)\\
							&\quad+ a([\xi_2,\varphi_2 Z_2]-\varphi_2^2 Z_2 -(a^2+b^2-1) \varphi_2^2 Z_2)\\
							& =0.
						\end{align*}
						Hence (i) is proved.
						
						For (ii), first note that $R(U_1,V_1)Z_1=R^1(U_1,V_1)Z_1$ follows from the fact that $\nabla_{X_1} Y_1=\nabla^1_{X_1} Y_1$ for $X_1,Y_1\in \X(S_1)$. For $Z_2$, we compute using Corollary \ref{corollary:nabla}
						
						\begin{align*}
							R(U_1,V_1)Z_2&=\nabla_{U_1} \nabla_{V_1} Z_2-\nabla_{V_1} \nabla_{U_1} Z_2-\nabla_{[U_1,V_1]} Z_2\\
							&=\nabla_{U_1} (-a \eta_2(Z_2)\varphi_1 V_1)-\nabla_{V_1}(-a\eta_2(Z_2)\varphi_1 U_1)\\
							&\quad+a[\eta_2(Z_2)\varphi_1[U_1,V_1]+\eta_1([U_1,V_1])\varphi_2 Z_2] \\
							&=a\eta_2(Z_2)(-\nabla_{U_1} \varphi_1 V_1+\nabla_{V_1} \varphi_1 U_1+\varphi_1[U_1,V_1])+a \eta_1([U_1,V_1]) \varphi_2 Z_2\\
							&=-a\eta_2(Z_2)(\Phi_1(U_1,\varphi_1 V_1)\xi_1-\nabla^{1,T}_{U_1} \varphi_1 V_1-\Phi(V_1,\varphi_1 U_1)\xi_1\\
							&\quad +\nabla^{1,T}_{V_1} \varphi_1 U_1 +\varphi_1 (\nabla^{1,T}_{U_1} V_1-\nabla^{1,T}_{V_1} U_1))-2a\Phi_1(U_1,V_1)\varphi_2 Z_2\\
							&=-2a \Phi_1(U_1,V_1)\varphi_2 Z_2,
						\end{align*}
						due to \eqref{eq:transverse-1}, where here $\nabla^{1,T}$ denotes the transverse connection on the Sasakian manifold $S_1$.
						
						For (iii), the computation of $R(U_2,V_2)Z_1$ is completely analogous to the computation of $R(U_1,V_1)Z_2$ in (ii) and so we omit it. 
						Finally, we compute 
						\begin{align*}
							R(U_2,V_2)Z_2  &= \nabla_{U_2}\nabla_{V_2}Z_2-\nabla_{V_2}\nabla_{U_2}Z_2 -\nabla_{[U_2,V_2]}Z_2 \\
							& = \nabla_{U_2}(\nabla^2_{V_2}Z_2-\lambda_{a,b}\eta_2(Z_2)\varphi_2 V_2) -\nabla_{V_2}(\nabla^2_{U_2}Z_2-\lambda_{a,b}\eta_2(Z_2)\varphi_2 U_2) \\
							& \quad -\nabla^2_{[U_2,V_2]}Z_2+\lambda_{a,b} [\eta_2(Z_2)\varphi_2[U_2,V_2]+\eta_2([U_2,V_2])\varphi_2 Z_2] \\
							& = \nabla^2_{U_2}\nabla^2_{V_2}Z_2-\lambda_{a,b}\eta_2(\nabla^2_{V_2}Z_2)\varphi_2U_2  \\
							&\quad -\lambda_{a,b} \left( U_2(\eta_2(Z_2))\varphi_2 V_2+\eta_2(Z_2)\nabla_{U_2}\varphi_2 V_2\right) \\
							& \quad -\nabla^2_{V_2}\nabla^2_{U_2}Z_2+\lambda_{a,b}\eta_2(\nabla^2_{U_2}Z_2)\varphi_2V_2  \\
							&\quad +\lambda_{a,b} \left( V_2(\eta_2(Z_2))\varphi_2 U_2 +\eta_2(Z_2)\nabla_{V_2}\varphi_2 U_2\right) \\
							& \quad -\nabla^2_{[U_2,V_2]}Z_2+\lambda_{a,b}[\eta_2(Z_2)\varphi_2[U_2,V_2]-2\Phi_2(U_2,V_2)\varphi_2 Z_2] \\
							& = R^2(U_2,V_2)Z_2-\lambda_{a,b}\left(g_2(\varphi_2 V_2,Z_2)\varphi_2 U_2 -g_2(\varphi_2 U_2,Z_2)\varphi_2V_2 \right)\\ 
							& \quad -\lambda_{a,b}\eta_2(Z_2)\left(-\Phi_2(U_2,\varphi_2 V_2)\xi_2 +\nabla^{2,T}_{U_2}\varphi_2 V_2 \right. \\ 
							& \quad \left. +\Phi_2(V_2,\varphi_2 U_2)\xi_2 -\nabla^{2,T}_{V_2}\varphi_2 U_2-\varphi_2(\nabla^{2,T}_{U_2}V_2-\nabla^{2,T}_{V_2}U_2) \right)\\
							& \quad -2\lambda_{a,b}\Phi_2(U_2,V_2)\varphi_2 Z_2 \\
							& = R^2(U_2,V_2)Z_2
							\\&\quad +\lambda_{a,b}(\Phi_2(V_2,Z_2)\varphi_2 U_2-\Phi_2(U_2,Z_2)\varphi_2 V_2-2\Phi_2(U_2,V_2)\varphi_2 Z_2).
						\end{align*}
						The last statement follows easily from (ii), (iii) and Lemma \ref{lemma:transverse}(iv).
					\end{proof}
					
					\medskip
					
					\section{Harmonicity of the complex structure \texorpdfstring{$J_{a,b}$}{} with respect to \texorpdfstring{$g_{a,b}$}{}}\label{sec:harmonic}
					
					Let $(M,g)$ denote a Riemannian manifold. According to \cite{Wood}, a $g$-orthogonal almost complex structure is harmonic if and only if 
					\begin{equation*}
						[J,\nabla^*\nabla J]=0, 
					\end{equation*}
					where $\nabla^*\nabla J$ is the \textit{rough Laplacian} of $J$ defined by $\nabla^*\nabla J=\Tr \nabla^2J$. That is, if $\{u_1,\ldots,u_{2n}\}$ is a local orthonormal frame on $M$, then 
					\[ (\nabla^*\nabla J)(W)= \sum_{i=1}^{2n} (\nabla^2_{u_i,u_i} J)(W), \quad W\in\mathfrak{X}(M), \] 
					where the second covariant derivative of $J$ is given by 
					\begin{equation*}
						(\nabla^2_{U,V} J)(W)=(\nabla_{U} (\nabla_{V} J))(W)-(\nabla_{\nabla_{U} V} J)(W).
					\end{equation*}
					It is clear that $\nabla^*\nabla J$ is a $(1,1)$-tensor on $M$.
					
					\medskip
					
					Let $(M^{2n},J,g)$ be an almost Hermitian manifold. According to \cite{Wood}, the following $2$-form $\rho$ plays a special role when determining if the almost complex structure $J$ is harmonic:
					\[ \rho=\mathcal{R}(\omega)\in\Omega^2(M),\]
					where $\omega$ is the fundamental $2$-form associated to $(J,g)$ and $\mathcal R$ is the curvature operator acting on $2$-forms. This $2$-form $\rho$ is a natural generalization of the Chern-Ricci form of a K\"ahler manifold,
					although in general it is not closed. It can be seen that the skew-symmetric tensor $P:TM\to TM$ obtained by contracting $\rho$ and $g$, i.e. $\rho(X,Y)=g(PX,Y)$, is given by \begin{equation}\label{eq:P}
						P(X)=\frac12 \sum_{i=1}^{2n} R(e_i, Je_i)X,\quad X\in\X(M),
					\end{equation}
					where $\{e_i\}$ is any orthonormal local frame of $M$. Also, let $\delta J\in\mathfrak{X}(M)$ denote the codifferential of $J$, that is, the unique vector field on $M$ satisfying 
					\[g(\delta J,X)= \delta\omega(X) \qquad \text{for all}\; X\in\X(M),\]
					where $\delta \omega$ is the codifferential of $\omega$. Since $\delta \omega$ is given by 
					\[ \delta \omega(X)=-\sum_{i=1}^{2n} (\nabla_{e_i} \omega)(e_i,X),\]
					for any local orthonormal frame $\{e_i\}$ of $M$, we obtain the following expression for $\delta J$:
					\begin{equation}\label{eq:delta-J}
						\delta J=\sum_{i=1}^{2n} (\nabla_{e_i} J)(e_i).
					\end{equation}
					With all these ingredients we may recall the following result from \cite{Wood}. In the Appendix we give an elementary proof of this fact.\footnote{The sign in the formula is different from \cite{Wood} since there $\omega=g(J\cdot, \cdot)$ but for us $\omega=g(\cdot,J\cdot)$.}

					\begin{proposition}\cite[Theorem 2.8]{Wood}\label{proposition:integrable-wood}
						Let $J$ be the almost complex structure of a $2n$-dimensional almost Hermitian manifold $(M,J,g)$. If $J$ is integrable then \[ [J,\nabla^*\nabla J]=2(\nabla_{\delta J} J-[J,P]).\]
						In particular, $J$ is harmonic if and only if $[J,P]=\nabla_{\delta J}J$.
					\end{proposition}
					
					The main result of this section is the following.
					
					\begin{theorem}\label{theorem:prod-sasakian}
						Let $(S_1^{2n_1+1},\varphi_1,\xi_1,\eta_1,g_1)$ and $(S_2^{2n_2+1},\varphi_2,\xi_2,\eta_2,g_2)$ be two Sasakian manifolds. If $(J,g):=(J_{a,b},g_{a,b})$ denotes the complex structure on $M_{a,b}=S_1\times S_2$ given in \eqref{eq:Jab} and \eqref{eq:gab} then $J$ is harmonic with respect to $g$, for any $a,b\in\R,\, b\neq 0$. 
					\end{theorem}
					
					\medskip
					
					First we prove an auxiliary result, which generalizes \cite[Lemma 5.4]{Tsu}, where the case of Calabi-Eckmann manifolds is considered.
					
					\begin{lemma}\label{lemma:codif-J}
						With the notation of Theorem \ref{theorem:prod-sasakian}, the codifferential $\delta J$ of the complex structure $J$ on $M_{a,b}$ is given by:
						\[ \delta J=2n_1\xi_1+2n_2\xi_2. \]
						Moreover, $\nabla_{\delta J}J=0$.
					\end{lemma}
					
					\begin{proof}
						We will compute $\delta J$ using \eqref{eq:delta-J}. Let us consider a local orthonormal frame on $M_{a,b}$ of the following form:
						\[ \left\{ \xi_1, J\xi_1=-\frac{a}{b}\xi_1+\frac{1}{b}\xi_2, e_1,\ldots,e_{2n_1}, f_1,\ldots, f_{2n_2} \right\}, \]
						where each $e_j$ is a local section of $\D_1$  and each $f_k$ is a local section of $\D_2$. With this frame, \eqref{eq:delta-J} becomes
						\begin{equation*}
							\delta J =(\nabla_{\xi_1} J)\xi_1+(\nabla_{J\xi_1} J)J\xi_1 +\sum_{j=1}^{2n_1} (\nabla_{e_j} J)e_j+\sum_{k=1}^{2n_2} (\nabla_{f_k} J)f_k.
						\end{equation*}
						Since $J$ is integrable, it follows from \cite[Corollary 2.2]{Gray} that $(\nabla_{J\xi_1} J)J\xi_1=(\nabla_{\xi_1} J)\xi_1$. Therefore, it follows from Lemma \ref{lemma:nabla-XX} that $(\nabla_{\xi_1} J)\xi_1=(\nabla_{J\xi_1} J)J\xi_1=0$, $(\nabla_{e_j} J)e_j=\xi_1$ and $(\nabla_{f_k} J)f_k=\xi_2$, so that
						\[ \delta J = \sum_{j=1}^{2n_1} \xi_1 +\sum_{k=1}^{2n_2} \xi_2= 2n_1\xi_1+2n_2\xi_2, \]
						as stated. 
						
						Finally, the last statement follows from Lemma \ref{lemma:nabla-XX}.
					\end{proof}
					
					\smallskip
					
					\begin{proof}[Proof of Theorem \ref{theorem:prod-sasakian}]
						As in Lemma \ref{lemma:codif-J}, we consider a local orthonormal frame on $M_{a,b}$ of the following form:
						\begin{equation}\label{eq:frame}
							\left\{ \xi_1, J\xi_1=-\frac{a}{b}\xi_1+\frac{1}{b}\xi_2, e_1,\ldots,e_{2n_1}, f_1,\ldots, f_{2n_2} \right\}, 
						\end{equation}
						where each $e_j$ is a local section of $\D_1$  and each $f_k$ is a local section of $\D_2$.
						
						Since $J$ is integrable, it follows from Proposition \ref{proposition:integrable-wood} and Lemma \ref{lemma:codif-J} that $J$ is harmonic if and only if $[J,P]=0$, with $P$ defined as in \eqref{eq:P} for the local frame \eqref{eq:frame}. We prove next that $J$ and $P$ commute indeed. We begin by computing $P$:
						\begin{align*}
							-2P & = 2R(\xi_1,J\xi_1)+\sum_j R(e_j,\varphi_1 e_j)+\sum_k R(f_k,\varphi_2 f_k) \\
							& =\sum_j R(e_j,\varphi_1 e_j)+\sum_k R(f_k,\varphi_2 f_k),
						\end{align*}
						due to Lemma \ref{lemma:curvatura}.
						
						Next we show that $[J,R(e_j,\varphi_1 e_j)]=[J,R(f_k,\varphi_2 f_k)]=0$. Recalling that $R(e_j,\varphi_1 e_j)\xi_i=R(f_k,\varphi_1 f_k)\xi_i=0$, $i=1,2$, it is enough to show that $[J,R(e_j,\varphi_1 e_j)]$ and $[J,R(f_k,\varphi_2 f_k)]$ vanish when evaluated in sections of $\D_1$ and $\D_2$. We use Lemma \ref{lemma:curvatura} for all the computations below. First, for $U_1\in \Gamma(\D_1)$, we compute
						\[ [J,R(e_j,\varphi_1 e_j)](U_1)=J(R^1(e_j,\varphi_1 e_j)U_1)-R^1(e_j,\varphi_1 e_j)\varphi_1 U_1=0,\]
						due to Corollary \ref{corollary:R-commutes}.
						Now, for $U_2\in\Gamma(\D_2)$,
						\begin{align*} 
							[J,R(e_j,\varphi_1 e_j)](U_2)&=2a\Phi_1(e_j,\varphi_1 e_j) U_2-R(e_j,\varphi_1 e_j) \varphi_2 U_2\\
							&=2a \Phi_1(e_j,\varphi_1 e_j)U_2-2a \Phi_1(e_j,\varphi_1 e_j) U_2\\
							&=0.
						\end{align*}
						Similarly, 
						\[[J,R(f_k,\varphi_2 f_k)](U_1)=2a \Phi_2(f_k,\varphi_2 f_k) U_1-2a \Phi_2(f_k,\varphi_2 f_k) U_1=0.\]
						Finally, 
						\begin{align*}
							[J,R(f_k,\varphi_2 f_k)](U_2)&=J(R(f_k,\varphi_2 f_k)U_2)-R(f_k,\varphi_2 f_k) \varphi_2 U_2\\
							&=J(R^2(f_k,\varphi_2 f_k)U_2)-R^2(f_k,\varphi_2 f_k)\varphi_2 U_2\\
							&\quad +\lambda_{a,b}(-\Phi_2(\varphi_2 f_k,U_2) f_k+\Phi_2(f_k,U_2) \varphi_2 f_k \\
							&\qquad\qquad +2\Phi_2(f_k,\varphi_2 f_k) U_2)  \\
							&\quad -\lambda_{a,b}(\Phi_2(\varphi_2 f_k, \varphi_2 U_2)\varphi_2 f_k+\Phi_2(f_k,\varphi_2 U_2) f_k\\
							&\qquad\qquad +2\Phi_2(f_k,\varphi_2 f_k)U_2)\\
							&=0,
						\end{align*}
						since $J(R^2(f_k,\varphi_2 f_k)U_2)=\varphi_2(R^2(f_k,\varphi_2 f_k)U_2)$ due to Corollary \ref{corollary:R-commutes}. It follows that $[J,P]=0$ and thus the complex structure $J=J_{a,b}$ is harmonic with respect to the metric $g=g_{a,b}$.
					\end{proof}
					
					\begin{remark}
						As mentioned before, Wood proved that Calabi-Eckmann manifolds equipped with the product of round metrics are harmonic. Moreover, he showed that, with the possible exception of $\mathbb{S}^1\times \mathbb{S}^3$ and $\mathbb{S}^3\times \mathbb{S}^3$, Calabi-Eckmann manifolds are \textit{unstable} (see \cite[\S7]{Wood}), that is, they do not correspond to a local minimum of the variational problem (or equivalently, if the second variation of the energy is positive, see \cite{Wood-Crelle}). However, his proof depends strongly on geometric properties of the spheres with the round metrics. This could indicate that the study of stability in general on a product of Sasakian manifolds is more complicated and it will be pursued in a future paper. 
					\end{remark}
					
					\medskip
					
					\section{Further properties of the Hermitian structures \texorpdfstring{$(J_{a,b},g_{a,b})$}{}}\label{sec:Hermitian}
					
					In this section we will focus on products of Sasakian manifolds which admit special Hermitian structures.
					
					Let $(M,J,g)$ be a Hermitian manifold with $\dim_\C M=n$. The complex structure $J$ can be extended naturally to differential forms on $M$ as follows: for a $p$-form $\alpha$, the $p$-form $J\alpha$ is given by
					\begin{align*}
						J\alpha=\alpha,&\quad p=0\\
						(J\alpha)(\cdot,\ldots,\cdot)=\alpha(J\cdot,\ldots,J\cdot),&\quad p>0.
					\end{align*}
					The real differential operator $d^c$ is then defined by 
					\begin{equation*}
						d^c \alpha=-J^{-1} d J\alpha=(-1)^p J d J\alpha, \quad \alpha\in \Omega^p(M).
					\end{equation*}
					It is well known that $dd^c=2\sqrt{-1}\partial \bar{\partial}$. 
					
					If the fundamental 2-form $\omega(\cdot,\cdot)=g(\cdot,J\cdot)$ is closed, then the metric $g$ is Kähler. The topology of compact K\"ahler manifolds is well understood, and there are several topological obstructions for the existence of a K\"ahler metric. 
					In the literature other conditions on the fundamental 2-form which are weaker than being closed have been introduced. 
					We recall next some of these non-K\"ahler Hermitian conditions. If $(M,J,g)$ is a Hermitian manifold with $\dim_\C M=n$, then the metric $g$ is said to be:
					\begin{enumerate}
						\item \textit{Balanced} if its fundamental form $\omega$ satisfies $d\omega^{n-1}=0$, or equivalently $\delta\omega=0$, where $\delta$ is the codifferential associated to $g$.
						\item \textit{Locally conformally Kähler} (LCK) if there exists an open covering $\{U_i\}_i$ of $M$ and smooth functions $f_i$ on each $U_i$ such that $e^{-f_i} g$ is Kähler. This definition is equivalent to the existence of a closed $1$-form $\theta$ such that $d\omega= \theta\wedge\omega$. The $1$-form $\theta$ coincides with the \textit{Lee form} associated to $(J,g)$, which is defined (on any Hermitian manifold) by: 
						\begin{equation}\label{eq:Lee}
							\theta=\frac{1}{n-1} (\delta\omega) \circ J.
						\end{equation} 
						If the Lee form of an LCK structure is parallel with respect to the Levi-Civita connection of $g$, the LCK metric is called \textit{Vaisman}. 
						\item \textit{Strong Kähler with torsion} (SKT), also called \textit{pluriclosed}, if its fundamental form satisfies $\partial\bar{\partial}\omega=0$, or equivalently, $dd^c \omega=0$. SKT metrics have applications in type II string theory and in 2-dimensional supersymmetric $\sigma$-models \cite{GHR,IP,St}. Moreover, they have also relations with generalized K\"ahler geometry (see for instance \cite{AG,Gua}).
						\item \textit{Astheno-Kähler} if its fundamental form satisfies $dd^c \omega^{n-2}=0$. Clearly this notion only makes sense for $n\geq 3$. Jost and Yau introduced these metrics in \cite{JY} to study Hermitian harmonic maps and to
						extend Siu’s rigidity theorem to non-K\"ahler manifolds.
						\item \textit{Gauduchon} if its fundamental form satisfies $dd^c \omega^{n-1}=0$. Any compact Hermitian manifold admits a Gauduchon metric in its conformal class, and it is unique in this class up to homotheties, due to a renowned result by Gauduchon \cite{Gau}.  
						\item $k$-\textit{Gauduchon}, for $1\leq k\leq n-1$, if its fundamental form $\omega$ satisfies 
						\[ dd^c \omega^k \wedge \omega^{n-k-1}=0.\] 
						These metrics were introduced in \cite{FWW} and they generalize the notion of Gauduchon metric, which corresponds precisely to $k=n-1$. Moreover, SKT and astheno-Kähler metrics are respectively 1-Gauduchon and $(n-2)$-Gauduchon. 
					\end{enumerate}
					
					\begin{remark}
						Notice that when $n=3$ a Hermitian metric is SKT if and only if it is astheno-K\"ahler, and in this case it is also $1$-Gauduchon.
					\end{remark}
					
					\medskip
					
					As in previous sections, $S_1$ and $S_2$ will be Sasakian manifolds with $\dim S_i=2n_i+1$, $n_i\in \N_0$, $i=1,2$, and the product manifold $M_{a,b}=S_1\times S_2$ will be equipped with the Hermitian structure $(J_{a,b},g_{a,b})$ defined in \eqref{eq:Jab} and \eqref{eq:gab}. 
					
					The fundamental $2$-form $\omega_{a,b}=g_{a,b}(\cdot, J_{a,b}\cdot )$ associated to  $(J_{a,b},g_{a,b})$
					is given by 
					\begin{equation}\label{eq:omega}
						\omega_{a,b}=\Phi_1+\Phi_2-b\eta_1\wedge \eta_2,
					\end{equation}
					where $\Phi_i$ and $\eta_i$ have been extended to the product in a natural way. Since both manifolds are Sasakian, the exterior derivative of $\omega_{a,b}$ is given by
					\begin{equation}\label{eq:d-omega}
						d\omega_{a,b}=-2b(\Phi_1\wedge \eta_2-\eta_1\wedge \Phi_2).
					\end{equation}
					
					Since $b\neq 0$, it is clear that $d\omega_{a,b}=0$ if and only if $\Phi_1=\Phi_2=0$, which can only happen when $n_1=n_2=0$, that is, $\dim S_1=\dim S_2=1$, hence $\dim M_{a,b}=2$. As we are interested in the non-Kähler setting, from now on we will assume that $n_1+n_2\geq 1$, so that $\dim M_{a,b}\geq 4$.
					
					Regarding the Hermitian structures $(J_{a,b},g_{a,b})$, Matsuo proved in \cite{Mat} that such a Hermitian structure is astheno-Kähler if and only if the following condition holds:
					\begin{equation}\label{eq:matsuo}
						n_1(n_1-1)+2a n_1n_2+n_2(n_2-1)(a^2+b^2)=0, 
					\end{equation} 
					provided that $n_1+n_2 \geq 2$, that is, $\dim_\C (S_1\times S_2)\geq 3$. Later, in \cite{FU}, Fino and Ugarte studied when the Hermitian structures $(J_{a,b},g_{a,b})$ are 1-Gauduchon and they obtained that this happens if and only if \eqref{eq:matsuo} holds, that is, if and only if they are astheno-Kähler. 
					
					\begin{remark}
						For any $n_1$ and $n_2$ such that $n_1+n_2\geq 2$ and $n_2\geq 1$ there are values of $a$ and $b$ such that $(J_{a,b},g_{a,b})$ is astheno-Kähler. Indeed, if $n_2=n_1=1$ then $a=0$, and if $n_2=1$, $n_1\geq 2$, then $a=-\frac{n_1-1}{2}$;  in both cases $b$ is arbitrary. Finally, if $n_2\geq 2$ the possible values of $a, b$ satisfy \[ a\in \left(-\frac{n_1}{n_2-1}-\frac{\sqrt{n_1n_2(n_1+n_2-1)}}{n_2(n_2-1)},-\frac{n_1}{n_2-1}+\frac{\sqrt{n_1n_2(n_1+n_2-1)}}{n_2(n_2-1)}  \right),\] 
						and 
						\[b^2=-a^2-\frac{2n_1n_2 a+n_1(n_1-1)}{n_2(n_2-1)}>0.\] 
						Note that if the Hermitian structure $(J_{a,b},g_{a,b})$ is astheno-Kähler then $a\leq 0$, and moreover, $a=0$ if and only if $n_1=n_2=1$. In particular, Morimoto's structure $(J_{0,1},g_{0,1})$ satisfies \eqref{eq:matsuo} if and only if $n_1=n_2=1$.
					\end{remark} 
					
					\medskip 
					
					In the next results we study when $(J_{a,b},g_{a,b})$ on $M_{a,b}=S_1\times S_2$ lies in one of the special classes of Hermitian structures mentioned above. 
					
					\begin{proposition}\label{proposition:balanced}
						The Hermitian structure $(J_{a,b},g_{a,b})$ on $M_{a,b}=S_1\times S_2$ with $n_1+n_2\geq 1$ is not balanced.
					\end{proposition}
					
					\begin{proof}
						The Hermitian structure is balanced if and only if $\delta \omega_{a,b}=0$, where $\omega_{a,b}$ is the associated fundamental $2$-form and $\delta$ is the codifferential. Equivalently, $\delta J_{a,b}=0$, where $\delta J_{a,b}$ is the vector field on $M_{a,b}$ dual to $\delta \omega_{a,b}$. However, it follows from Lemma \ref{lemma:codif-J} that $\delta J_{a,b}=2n_1\xi_1+2n_2\xi_2\neq 0$ since $n_1+n_2\geq 1$.
					\end{proof}
					
					\begin{remark}
						Some products $S_1\times S_{2}$ of Sasakian manifolds carry balanced Hermitian structures, different from any $(J_{a,b},g_{a,b})$. For instance, if $H_3$ is the $3$-dimensional Heisenberg group, it is well known that $H_3$ carries a natural left invariant Sasakian structure (see Table \ref{table:Sasakian 3 Lie algebras} in \S\ref{section:examples} below). On the other hand, it was shown in \cite{FPS,Ug} that $G:=H_3\times H_3$ carries a left invariant balanced Hermitian structure. If $\Gamma_1$ and $\Gamma_2$ are co-compact discrete subgroups of $H_3$ then $\Gamma:=\Gamma_1\times\Gamma_2$ is a discrete co-compact subgroup of $G$. Therefore the nilmanifold $\Gamma\backslash G\cong (\Gamma_1\backslash H_3)\times (\Gamma_2\backslash H_3)$ is a product of Sasakian $3$-manifolds and carries a balanced Hermitian structure. 
					\end{remark}
					
					\medskip 
					
					\begin{proposition}\label{proposition:lck}
						The Hermitian structure $(J_{a,b},g_{a,b})$ on $M_{a,b}=S_1\times S_2$ is non-K\"ahler LCK if and only if $\dim S_1=1$ and $\dim S_2 \geq 3$ or $\dim S_2=1$ and $\dim S_1 \geq 3$. Moreover, the LCK structure is Vaisman.
					\end{proposition}
					
					\begin{proof}
						According to \eqref{eq:Lee}, the Lee form $\theta$ associated to $(J_{a,b},g_{a,b})$ is given by
						\[ \theta=\frac{1}{n_1+n_2}\delta \omega_{a,b} \circ J_{a,b}. \]
						Due to \ref{lemma:codif-J}, together with \eqref{eq:Jab} and \eqref{eq:gab}, we arrive easily at the following expression for $\theta$:
						\begin{equation}\label{eq:theta}
							\theta= \frac{2b}{n_1+n_2}(n_2\eta_1-n_1\eta_2). 
						\end{equation}
						Replacing \eqref{eq:omega}, \eqref{eq:d-omega} and \eqref{eq:theta}  in the LCK condition $d\omega=\theta\wedge \omega$,  we arrive at
						\[ -n_2\Phi_1\wedge \eta_2+n_1 \eta_1\wedge \Phi_2= n_2\eta_1\wedge \Phi_1-n_1\eta_2\wedge \Phi_2.\]
						Comparing the components in $S_1$ and $S_2$ of each term in this equation, we observe that all of them have to vanish, so that
						\[ n_2\Phi_1\wedge \eta_2 = n_1 \eta_1\wedge \Phi_2= n_2\eta_1\wedge \Phi_1 = n_1\eta_2\wedge \Phi_2=0. \]
						If $n_2 \neq 0$ then $\Phi_1\wedge \eta_2=\eta_1\wedge \Phi_1=0$ and this happens if and only if $\Phi_1=0$, or equivalently $\dim S_1=1$. In a similar way, if $n_1\neq 0$ then $\dim S_2=1$. 
						
						The fact the LCK structure is actually Vaisman, that is, $\nabla \theta=0$, follows immediately from Corollary \ref{corollary:nabla}.
					\end{proof}
					
					\medskip 
					
					Next we will analyze the Hermitian conditions which involve the operator $d^c$. First, it was shown in \cite{Mat} that $J_{a,b}\Phi_1=\Phi_1$ and $J_{a,b}\Phi_2=\Phi_2$. On the other hand, it is easy to verify that 
					\[ J_{a,b}\eta_1=\frac{a}{b}\eta_1+\frac{a^2+b^2}{b}\eta_2, \quad  J_{a,b}\eta_2=-\frac{1}{b}\eta_1-\frac{a}{b}\eta_2.\] 
					Using this, together with \eqref{eq:d-omega} and the definition of $d^c$, the following expressions are obtained:
					\begin{align}
						d^c \omega_{a,b}&=2[\Phi_1 \wedge (\eta_1+a\eta_2)+\Phi_2 \wedge (a\eta_1+(a^2+b^2)\eta_2)], \nonumber\\
						dd^c \omega_{a,b}&=4[\Phi_1^2+2a \Phi_1 \wedge \Phi_2+(a^2+b^2)\Phi_2^2], \label{eqs:dc-omega}\\
						d\omega_{a,b} \wedge d^c \omega_{a,b}&=4b[\Phi_1^2+2a \Phi_1\wedge \Phi_2 +(a^2+b^2)\Phi_2^2] \wedge \eta_1 \wedge \eta_2 \nonumber \\
						& = b\, dd^c\omega_{a,b} \wedge \eta_1 \wedge \eta_2.\nonumber
					\end{align}
					
					\begin{proposition}\label{proposition:skt}
						The Hermitian structure $(J_{a,b},g_{a,b})$ on $M_{a,b}=S_1\times S_2$ is non-K\"ahler SKT if and only if:
						\begin{enumerate} 
							\item [$(a)$] $\dim M_{a,b}=4$, or 
							\item [$(b)$] $\dim S_1=\dim S_2=3$ and $a=0$. 
						\end{enumerate}
					\end{proposition}
					
					\begin{proof}
						The Hermitian structure is SKT if and only if $dd^c \omega_{a,b}=0$. It follows from \eqref{eqs:dc-omega} that this happens if and only if \[\Phi_1^2+2a \Phi_1 \wedge \Phi_2+(a^2+b^2)\Phi_2^2=0.\]
						Comparing the components in $S_1$ and $S_2$ of each term in this equation, we deduce that the SKT condition is equivalent to  
						\begin{equation}\label{eqs:skt} 
							\Phi_1^2=0, \quad \Phi_2^2=0 \quad \text{and}\quad a \Phi_1 \wedge \Phi_2=0.
						\end{equation}
						If \eqref{eqs:skt} holds, from the first two equalities we obtain that $\dim S_i\leq 3$ for $i=1,2$. Since the metric is non-Kähler, $\dim S_i=3$ for $i=1$ or $i=2$. If both $\dim S_1=\dim S_2=3$, it follows from $a \Phi_1 \wedge \Phi_2=0$ that $a=0$.
						
						Conversely, if the conditions in the statement hold then it is clear that \eqref{eqs:skt} is satisfied.
					\end{proof}
					
					\begin{remark}
						It follows from Proposition \ref{proposition:skt} that if $S_1\times S_2$ admits a non-Kähler SKT structure then at least one of the factors has dimension 3. According to \cite{Geiges}, any compact $3$-dimensional Sasakian manifold is diffeomorphic to one of the following manifolds:
						\[ \mathbb{S}^3/\Gamma, \qquad H_3/\Gamma, \qquad \widetilde{\operatorname{SL}}(2,\R)/\Gamma,   \]
						where $\Gamma$ is a discrete subgroup of the isometry group of the corresponding canonical Sasakian metric. 
					\end{remark}

					\begin{remark}
						Note that according to Proposition \ref{proposition:lck} and Proposition \ref{proposition:skt}, the Hermitian structures $(J_{a,b}, g_{a,b})$ on the complex surfaces $S_1\times S_2$ with $\dim S_1=1$ and $\dim S_2=3$ are both Vaisman and SKT. This is not a surprise since it was proved in \cite[Theorem A]{FT} that a Hermitian metric on a complex surface is Vaisman if and only if the metric is SKT and the Bismut connection satisfies the first Bianchi identity. The fact that the Bismut connection associated to $(J_{a,b},g_{a,b})$ on the complex surface $S_1\times S_2$ satisfies the first Bianchi identity follows from \cite[Theorem 3.2]{FT} and \cite[Proposition 3.2]{Bel}.
					\end{remark}
					
					\medskip 
					
					In the following theorem we characterize when the Hermitian structure $(J_{a,b},g_{a,b})$ is $k$-Gauduchon for $\dim_\C M_{a,b}\geq 4$. 
					Since the case $k=1$ was already done in \cite{FU}, we will restrict to the case $2\leq k \leq n-1$. We obtain that the metric $g_{a,b}$ is always Gauduchon and, furthermore, it is $k$-Gauduchon with $2\leq k\leq n-1$ if and only if it is astheno-Kähler, just as in the case $k=1$.
					
					\begin{theorem}\label{theorem:k-gau}
						Let $M_{a,b}=S_1\times S_2$ be a product of Sasakian manifolds with $n:=n_1+n_2+1\geq 4$. Then the Hermitian structure $(J_{a,b},g_{a,b})$ on $M_{a,b}$ is $k$-Gauduchon with $2\leq k\leq n-1$ if and only if the constants $a$ and $b$ satisfy 
						\begin{equation}\label{eq:gauduchon}
							(n-1-k) [n_1(n_1-1)+2a n_1n_2+n_2(n_2-1)(a^2+b^2)]=0.
						\end{equation}
						In particular, $(J_{a,b}, g_{a,b})$ is Gauduchon and it is $k$-Gauduchon with $2\leq k\leq n-2$ if and only if it is astheno-Kähler.
					\end{theorem}
					
					\begin{proof}
						Let us denote $J=J_{a,b}$ and $\omega=\omega_{a,b}$. We follow the lines of the proof of \cite[Theorem 4.1]{Mat}. 
						We will use equations \eqref{eqs:dc-omega} and the fact that $J\omega^k=\omega^k$ for all $k$ since $\omega$ is a $(1,1)$-form on $M$. For $2\leq k\leq n-1$, we compute
						\begin{align*}
							dd^c \omega^k \wedge \omega^{n-k-1}&=d(JdJ\omega ^k) \wedge \omega^{n-k-1}\\
							&= k\,d(J( \omega^{k-1} \wedge d\omega)) \wedge \omega^{n-k-1}\\
							&=k\,d(\omega^{k-1} \wedge d^c \omega) \wedge \omega^{n-k-1}  \\ &=k [(k-1) \omega^{k-2} \wedge d\omega \wedge d^c \omega+\omega^{k-1} \wedge dd^c  \omega] \wedge \omega^{n-k-1} \\
							&=k [(k-1) d\omega \wedge d^c\omega+\omega \wedge dd^c\omega] \wedge \omega^{n-3}\\
							&=k \, dd^c\omega \wedge [b(k-2) \eta_1\wedge \eta_2+\Phi_1+\Phi_2] \wedge \omega^{n-3}.
						\end{align*}
						Since $(\eta_1 \wedge \eta_2)^2=0$, it follows from the binomial theorem that
						\begin{align*}
							\omega^{n-3} & = (\Phi_1+\Phi_2-b\, \eta_1\wedge\eta_2)^{n-3} \\
							& = (\Phi_1+\Phi_2)^{n-3}-(n-3)(\Phi_1+\Phi_2)^{n-4} \wedge (b\, \eta_1 \wedge \eta_2),
						\end{align*}
						and therefore 
						\begin{align*}
							dd^c \omega^k \wedge \omega^{n-k-1}&=k\, dd^c\omega \wedge [b(k-n+1)(\Phi_1+\Phi_2)^{n-3} \wedge \eta_1 \wedge \eta_2+(\Phi_1+\Phi_2)^{n-2}]\\
							&=k\, dd^c\omega \wedge [b(k-n+1) \eta_1 \wedge \eta_2+(\Phi_1+\Phi_2)] \wedge (\Phi_1+\Phi_2)^{n-3}.
						\end{align*}
						Given that $\Phi_1^p=0$ when $p>n_1$ and $\Phi_2^p=0$ when $p>n_2$, an index $j$ satisfies $j\leq n_1$ and $n-3-j \leq n_2$ only when $n_1-2 \leq j \leq n_1$. Hence,
						\[ (\Phi_1+\Phi_2)^{n-3}=\binom{n-3}{n_1-2} \Phi_1^{n_1-2} \wedge \Phi_2^{n_2}+\binom{n-3}{n_1-1} \Phi_1^{n_1-1} \wedge \Phi_2^{n_2-1}+\binom{n-3}{n_1} \Phi_1^{n_1} \wedge \Phi_2^{n_2-2}.\]
						Therefore, using the expression for $dd^c \omega$ given in \eqref{eqs:dc-omega},
						\[ dd^c \omega^k \wedge \omega^{n-k-1}=4bk(k-n+1) C(n,n_1) \Phi_1^{n_1}\wedge \Phi_2^{n_2}\wedge \eta_1 \wedge \eta_2,\] where $C(n,n_1)=\binom{n-3}{n_1-2}+2a \binom{n-3}{n_1-1}+(a^2+b^2)\binom{n-3}{n_1}$.
						
						Since $4bk\neq 0$ and $\Phi_1^{n_1}\wedge \Phi_2^{n_2}\wedge \eta_1 \wedge \eta_2$ is a volume form on $M_{a,b}$,  we get that $dd^c\omega^k\wedge\omega^{n-k-1}=0$  if and only if $(k-n+1) C(n,n_1)=0$, and this is equivalent to \eqref{eq:gauduchon}.
					\end{proof}
					
					\smallskip
					
					Let us now analyze the $k$-Gauduchon condition in the case $\dim_\C M_{a,b}=3$. In \cite{FU} it was proved that the metric $(J_{a,b},g_{a,b})$ is 1-Gauduchon if and only if it is SKT if and only if $a=0$. We deal next with the only missing case $k=2$, which corresponds to Gauduchon metrics.
					
					\begin{proposition}
						Let $M_{a,b}=S_1\times S_2$ be a product of Sasakian manifolds with $\dim_\C M_{a,b}=3$. Then the Hermitian structure $(J_{a,b},g_{a,b})$ is Gauduchon. 
					\end{proposition}
					
					\begin{proof} 
						It follows from the computations in the proof of Theorem \ref{theorem:k-gau} that 
						\begin{align*}
							dd^c \omega^2&= 2 \, dd^c \omega \wedge (\Phi_1+\Phi_2)\\
							&=8[  \Phi_1^3+(2a+1) \Phi_1^2 \wedge \Phi_2+(a^2+b^2+2a) \Phi_1 \wedge \Phi_2^2+(a^2+b^2) \Phi_2^3],
						\end{align*}
						using \eqref{eqs:dc-omega}.
						We have to analyze two cases: (i) $\dim S_1=\dim S_2=3$ and (ii) $\dim S_1=1$ and $\dim S_2=5$. For (i), note that $\Phi_i^2=0$ for $i=1,2$ and therefore $dd^c\omega^2=0$. For (ii), we have that $\Phi_1=0$ and $\Phi_2^3=0$ and therefore again $dd^c\omega^2=0$.
					\end{proof}
					
					\medskip
					
					\section{The Bismut connection on \texorpdfstring{$M_{a,b}=S_1\times S_2$}{}}\label{sec:Bismut}
					
					In this section we exhibit an explicit formula for the Bismut connection $\nabla^B$ on $M_{a,b}=S_1\times S_2$ associated to the Hermitian structure $(J_{a,b}, g_{a,b})$, in terms of the characteristic connections of the Sasakian manifolds $S_1$ and $S_2$.  As an application we study the Ricci curvature $\Ric^B$ and the Ricci form $\rho^B$ associated to $\nabla^B$. In particular, we characterize the Hermitian structures such that $\Ric^B=0$ and those known as Calabi-Yau with torsion, i.e. $\rho^B=0$, in terms of Sasakian $\eta$-Einstein manifolds. We use this characterization to provide examples of such manifolds.
					
					\smallskip 
					
					Let us recall some basic facts about metric connections with totally skew-symmetric torsion. For more details we refer to \cite{Ag}.
					
					A metric connection $D$ with torsion $T$ on a Riemannian manifold $(M,g)$ is said to have \emph{totally skew-symmetric torsion},
					or \emph{skew torsion} for short, if the $(0,3)$-tensor field $T$ defined by
					\[T(X,Y,Z)=g(T(X,Y),Z)\]
					is a $3$-form. The relation between $D$ and the Levi-Civita connection
					$\nabla$ is then given by
					\begin{equation}\label{eq:skew}
						D_XY=\nabla_XY+\frac{1}{2}T(X,Y).
					\end{equation}
					It is well known that $D$ has the same geodesics as $\nabla$. 
					
					\medskip 
					
					Both on Hermitian manifolds and Sasakian manifolds there is a distinguished metric connection with skew-symmetric torsion:
					\begin{enumerate}
						\item On any Hermitian manifold $(M,J,g)$ there exists a unique metric connection $\nabla^B$ with skew-symmetric torsion such that $\nabla^B J=0$ (see \cite{Bi}). This connection is known as the \textit{Bismut} connection (or \textit{Strominger}) and its torsion $3$-form is given by $T^B(X,Y,Z)=d^c\omega(X,Y,Z)=-d\omega(JX,JY,JZ)$, where $\omega=g(\cdot,J\cdot)$ is the fundamental $2$-form.
						\item On any Sasakian manifold $(S,\varphi,\xi,\eta,g)$ there exists a unique metric connection $\nabla^C$ with skew-symmetric torsion such that $\nabla^C\varphi=\nabla^C\eta=\nabla^C\xi=0$ (see \cite{FI}). It is known as the \textit{characteristic} connection\footnote{In fact, it was proved in \cite{FI} that the characteristic connection exists for a larger class of almost contact metric manifolds, namely, those which satisfy that $N_\varphi$ is totally skew-symmetric and $\xi $ is a Killing vector field.} and its torsion is given by $T^C=\eta\wedge d\eta$. If we consider $T^C$ as the usual $(1,2)$-tensor $T^C(X,Y)=\nabla^C_XY-\nabla^C_YX-[X,Y]$, we obtain that $T^C$ is given by
						\begin{equation}\label{eq:tc}
							T^C(X,Y)=2(-\eta(X)\varphi Y+\eta(Y)\varphi X+\Phi(X,Y)\xi).
						\end{equation}
					\end{enumerate}
					
					\medskip
					
					We compute now the Bismut connection on $M_{a,b}=S_1\times S_2$ associated to $(J_{a,b},g_{a,b})$.
					
					\begin{proposition}\label{proposition:bismut} The Bismut connection $\nabla^B$ associated to the Hermitian structure $(J_{a,b},g_{a,b})$ on the product $M_{a,b}=S_1\times S_2$ of two Sasakian manifolds is given by:
						\begin{enumerate}
							\item[$\ri$] $\nabla^B_{X_1} Y_1=\nabla^{1,C}_{X_1} Y_1\in \X(S_1)$,
							\item[$\rii$] $\nabla^B_{X_2}Y_2=\nabla^{2,C}_{X_2}Y_2-2(a^2+b^2-1) \eta_2(X_2)\varphi_2Y_2\in\X(S_2)$,
							\item[$\riii$] $\nabla^B_{X_1}Y_2=-2a \eta_1(X_1)\varphi_2Y_2\in \X(S_2)$,
							\item[$\riv$] $\nabla^B_{X_2}Y_1=-2a\eta_2(X_2)\varphi_1Y_1\in \X(S_1)$,
						\end{enumerate}
						where $\nabla^{i,C}$ is the characteristic connection on the Sasakian manifold $S_i$, $i=1,2$. 
						
						In particular, $\nabla^B \xi_1=\nabla^B\xi_2=0$. Moreover, the Lee form $\theta$ associated to $(J_{a,b},g_{a,b})$ is $\nabla^B$-parallel, $\nabla^B\theta=0$.
					\end{proposition}
					
					\begin{proof}
						From \eqref{eq:d-omega} we obtain the following formula for $T^B$:
						\begin{equation}\label{eq:T^B}
							T^B=2[\Phi_1 \wedge (\eta_1+a \eta_2)+\Phi_2 \wedge (a\eta_1+(a^2+b^2)\eta_2)].
						\end{equation}
						Replacing $T^B$ in \eqref{eq:skew} and using the expressions for the Levi-Civita connection obtained in Corollary \ref{corollary:nabla} we obtain $\ri$-$\riv$.
						
						It is immediate to verify from $\ri$-$\riv$ that both $\xi_1$ and $\xi_2$ are $\nabla^B$-parallel. From this, together with the fact that $\eta$ is $\nabla^C$-parallel on  any Sasakian manifold, we obtain that $\nabla^B \eta_1=\nabla^B \eta_2=0$. Then it follows from \eqref{eq:theta} that 
						$\nabla^B \theta=0$.
					\end{proof}
					
					\begin{corollary}
						On $M_{0,1}=S_1\times S_2$ equipped with the Morimoto's structure $(J_{0,1}, g_{0,1})$, the associated Bismut connection satisfies 
						\[ \nabla^B_{X_1+X_2} (Y_1+Y_2)=\nabla^{1,C}_{X_1} Y_1+\nabla^{2,C}_{X_2}Y_2.\]
					\end{corollary}

					\medskip
					
					In what follows we will study the vanishing of the Ricci curvatures associated to $\nabla^B$. Let us recall first that on a Hermitian manifold $(M^{2n},J,g)$ the Ricci tensor $\Ric^B$ and the Ricci form $\rho^B$ of the Bismut connection $\nabla^B$ are defined by
					\begin{equation}\label{eq:Ric-rho}
						\Ric^B(X,Y)=\sum_{i=1}^{2n} g(R^B(u_i,X)Y,u_i), \quad \rho^B(X,Y)=\frac{1}{2} \sum_{i=1}^{2n} g(R^B(X,Y)u_i,Ju_i),
					\end{equation} 
					where $\{u_1,\ldots, u_{2n}\}$ is a local orthonormal frame of $M$. We will compute $\Ric^B$ and $\rho^B$ using formulas appearing in \cite{IP}, and as a consequence we will have no need to determine explicitly the Bismut curvature tensor $R^B$.
					
					In this article we will use the following convention, as in \cite{GGP,Grant}. We will say that the Hermitian structure $(J,g)$ on $M$ is \textit{Calabi-Yau with torsion} (CYT, for short) if the Ricci form $\rho^B$ associated to the Bismut connection vanishes identically, or equivalently, the (restricted) holonomy group of the Bismut connection is contained in $\operatorname{SU}(n)$. 
					
					
					\medskip
					
					Resuming our study of the Bismut connection on the product of Sasakian manifolds, we recall the following result in \cite{Bel}, which can also be verified using Proposition \ref{proposition:bismut}.
					
					\begin{proposition}\cite[Proposition 3.2]{Bel}\label{proposition:torsion-paralela}
						The Bismut connection $\nabla^B$ associated to $(J_{a,b},g_{a,b})$ on $M_{a,b}=S_1\times S_2$ has parallel torsion, i.e., $\nabla^B T^B=0$.
					\end{proposition}
					
					As a consequence, it follows from \cite{CMS} that the Bismut curvature $R^B$ associated to  $(J_{a,b},g_{a,b})$ satisfies 
					\[ g(R^B(X,Y)Z,W)=g(R^B(Z,W)X,Y)\] 
					and from this it can be readily seen that the Bismut-Ricci tensor $\Ric^B$ is symmetric. According to \cite{IP} this is equivalent to 
					\begin{equation}\label{eq:codif T^B} 
						\delta T^B=0,
					\end{equation} where $\delta$ is the codifferential and $T^B$ is the torsion 3-form. 
					
					\medskip 
					
					\begin{remark}
						Proposition \ref{proposition:torsion-paralela} allows us to determine which of the Hermitian structures $(J_{a,b},g_{a,b})$ on $M_{a,b}$ are \textit{Kähler-like}, i.e. they satisfy, for any vector fields $X,Y,Z,W$, 
						\begin{itemize}
							\item (First Bianchi identity) $R^B(X,Y,Z)+R^B(Y,Z,X)+R^B(Z,X,Y)=0$,
							\item $R^B(X,Y,Z,W)=R^B(JX,JY,Z,W)$.
						\end{itemize}
						In fact, Proposition \ref{proposition:torsion-paralela} implies the second condition due to  \cite[Lemma 3.13]{AV1}. Moreover, due to \cite[Theorem 3.2]{FT}, the Bismut connection associated to $(J_{a,b},g_{a,b})$ on $M_{a,b}$ satisfies the first Bianchi identity if and only if $(J_{a,b},g_{a,b})$ is SKT. Therefore, $(J_{a,b},g_{a,b})$ is Kähler-like if and only if it is SKT.
					\end{remark}
					
					We proceed now to study when the Hermitian structure $(J_{a,b},g_{a,b})$ on $M_{a,b}=S_1\times S_2$ satisfies either $\Ric^B=0$ or $\rho^B=0$. We will show that these structures are closely related to a special family of Sasakian manifolds, namely, the $\eta$-Einstein ones.
					A $(2n+1)$-dimensional Sasakian manifold $(S,\varphi,\xi,\eta,g)$ is said to be \textit{$\eta$-Einstein} if the Ricci curvature tensor of
					the metric $g$ satisfies the equation 
					\begin{equation}\label{eq:eta-Einstein} 
						\Ric = \lambda g + \nu \eta \otimes \eta,
					\end{equation} 
					for some constants $\lambda, \nu\in \R$. 
					
					It is well known that on any Sasakian manifold $S$ the Riemannian curvature tensor $R$ satisfies \[R(X,Y)\xi=\eta(Y)X-\eta(X)Y,\]
					for all $X,Y\in\X(S)$. 
					It follows easily from this equation that $\Ric(\xi,X)=2n\eta(X)$ for any vector field $X$ on $S$ and using this fact we obtain:
					\begin{itemize}
						\item $\lambda + \nu=2n$,
						\item the $\eta$-Einstein condition \eqref{eq:eta-Einstein} with constants $(\lambda,2n-\lambda)$ is equivalent to 
						\begin{equation}\label{eq:Ric-D}
							\Ric(U,V)=\lambda g(U,V), \qquad U,V\in \Gamma(\mathcal D).
						\end{equation}
					\end{itemize} 
					Another immediate consequence of the definition is that every $\eta$-Einstein manifold is necessarily of constant scalar curvature $s=2n(\lambda+1)$.
					
					\smallskip
					
					In the next results we characterize the Hermitian structures $(J_{a,b},g_{a,b})$ such that $\Ric^B=0$. We will need an explicit expression for the (1,2)-tensor $T^B$, which we summarize in the following lemma, whose proof follows from Proposition \ref{proposition:bismut} and equation \eqref{eq:tc}.
					
					\begin{lemma}\label{lemma:T^B}
						Let $X_i\in \X(S_i)$ and $\{\xi_1,J\xi_1,e_1,\ldots,e_{2n_1},f_1,\ldots,f_{2n_2}\}$ be a local orthonormal frame on $M_{a,b}=S_1\times S_2$ as in \eqref{eq:frame}. Then,
						\begin{align*}
							T^B(X_1,\xi_1)&=2\varphi_1 X_1, \quad T^B(X_1,J\xi_1)=0\\
							T^B(X_1,e_j)&=-2(\eta_1(X_1)\varphi_1 e_j+\Phi_1(e_j,X_1)\xi_1),  \\ 
							T^B(X_1,f_k)&=-2a \eta_1(X_1)\varphi_2 f_k,\\
							T^B(X_2,\xi_1)&=2a \varphi_2 X_2, \quad T^B(X_2,J\xi_1)=2b  \varphi_2 X_2,  \\
							T^B(X_2,e_j)&=-2a \eta_2(X_2)\varphi_1 e_j,  \\
							T^B(X_2,f_k)&=-2(a^2+b^2)\eta_2(X_2)\varphi_2 f_k+2 \Phi_2(X_2,f_k)\xi_2. 
						\end{align*}
						In particular, for $X_i,Y_i\in \Gamma(\D_i)$ we obtain 
						\[ T^B(X_i,Y_i)=2\Phi_i(X_i,Y_i)\xi_i, \quad \text{and} \quad T^B(X_1,Y_2)=0.\]
					\end{lemma}
					
					\smallskip
					
					\begin{theorem}\label{theorem:BRF}
						The Bismut-Ricci tensor $\Ric^B$ associated to the Hermitian structure $(J_{a,b},g_{a,b})$ is given by:
						\begin{align*}
							\Ric^B(X_1,Y_1)&=\Ric^1(X_1,Y_1)-2g_1(X_1,Y_1)-(2n_1-2)\eta_1(X_1)\eta_1(Y_1)\\
							\Ric^B(X_1,Y_2)&=0\\
							\Ric^B(X_2,Y_2)& = \Ric^2(X_2,Y_2)-2(2a^2+2b^2-1)g_2(X_2,Y_2) \\
							& \quad +[2(2a^2+2b^2-1)-2n_2) ]\eta_2(X_2)\eta_2(Y_2),
						\end{align*} where $X_i,Y_i\in\X(S_i)$, for $i=1,2$ and $\Ric^i$ denotes the Ricci curvature associated to the Levi-Civita connection $\nabla^i$ on $S_i$. In particular, $\Ric^B=0$ if and only if both $S_1$ and $S_2$ are $\eta$-Einstein with constants 
						\begin{align*} 
							(\lambda_1,\nu_1)&=(2,2n_1-2),\\
							(\lambda_2,\nu_2)&=(2(2a^2+2b^2-1),2n_2-2(2a^2+2b^2-1)),
						\end{align*}
						respectively.
					\end{theorem}
					
					\begin{proof}
						
						To compute the Bismut-Ricci tensor $\Ric^B$ associated to $(J_{a,b},g_{a,b})$ we will use the formula for $\Ric^B$ on a general Hermitian manifold $(M^{2n},g)$ obtained in \cite[Proposition 3.1]{IP}:
						\[  \Ric(X,Y)=\Ric^B(X,Y)+\frac{1}{2} \delta T^B(X,Y)+\frac{1}{4} \sum_{i=1}^{2n} g(T^B(X,u_i),T^B(Y,u_i)),\] where $\Ric$ denotes the Ricci tensor of $g$ and $\{u_i\}_{i=1}^{2n}$ is a local orthonormal frame. 
						
						Using that $\delta T^B=0$ for the Hermitian structure $(J_{a,b},g_{a,b})$ on $M_{a,b}=S_1\times S_2$, in the local orthonormal frame \eqref{eq:frame} the expression of $\Ric^B$ simplifies to 
						\begin{align}\label{eq:Ric^B}
							\Ric^B(X,Y)&=\Ric(X,Y)-\frac14 \big[g(T^B(X,\xi_1),T^B(Y,\xi_1))+g(T^B(X,J\xi_1),T^B(Y,J\xi_1))\\
							&\quad +\sum_{j=1}^{2n_1} g(T^B(X,e_j),T^B(Y,e_j))+\sum_{k=1}^{2n_2} g(T^B(X,f_k),T^B(Y,f_k))\big], \nonumber
						\end{align}
						for $X,Y\in\X(M_{a,b})$.
						In \cite{LPS} the Ricci curvature of the metric $g_{a,b}$ on $M_{a,b}$ has been computed:
						\begin{align*}
							\Ric(X_1,Y_1)&=\Ric^1(X_1,Y_1)+2a^2 n_2 \eta_1(X_1)\eta_1(Y_1) \nonumber\\
							\Ric(X_1,Y_2)&=2a(n_1+n_2(a^2+b^2))\eta_1(X_1)\eta_2(Y_2)\\
							\Ric(X_2,Y_2)&=\Ric^2(X_2,Y_2)-2\lambda_{a,b}g_2(X_2,Y_2) +2[n_1 a^2+\lambda_{a,b}+n_2(a^2+b^2)^2-n_2]\eta_2(X_2)\eta_2(Y_2),\nonumber
						\end{align*}
						where $\lambda_{a,b}=a^2+b^2-1$.
						
						Next we compute $\Ric^B$ using \eqref{eq:Ric^B}, the formulas for $\Ric$ above and Lemma \ref{lemma:T^B}. For $X_1,Y_1\in \X(S_1)$, we have that
						\begin{align*}
							&\Ric^B(X_1,Y_1)=\Ric^1(X_1,Y_1)+2a^2 n_2 \eta_1(X_1)\eta_1(Y_1)-\frac14 \big[ 4 g_1(\varphi_1 X_1, \varphi_1 Y_1)\\
							&\quad +4\sum_{j=1}^{2n_1} g_1(\eta_1(X_1)\varphi_1 e_j+\Phi_1(e_j,X_1)\xi_1, \eta_1(Y_1)\varphi_1 e_j+\Phi_1(e_j,Y_1)\xi_1)\\
							&\quad + \sum_{k=1}^{2n_2} g_2(-2a \eta_1(X_1)\varphi_2 f_k, -2a \eta_1(Y_1) \varphi_2 f_k) \big]\\
							&=\Ric^1(X_1,Y_1)+2a^2 n_2 \eta_1(X_1)\eta_1(Y_1)-\frac14\big[ 4g_1(\varphi_1 X_1, \varphi_1 Y_1) +8 n_1 \eta_1(X_1)\eta_1(Y_1)\\
							&\quad +4\sum_{j=1}^{2n_1} g_1(e_j,\varphi_1 X_1)g_1(e_j,\varphi_1 Y_1) +8 n_2 a^2 \eta_1(X_1)\eta_1(Y_1) \big] \\
							&=\Ric^1(X_1,Y_1)-g_1(\varphi_1 X_1,\varphi_1 Y_1)-2 n_1 \eta_1(X_1)\eta_1(Y_1)-g_1(\varphi_1 X_1,\varphi_1 Y_1)\\
							&=\Ric^1(X_1,Y_1)-2g_1(X_1,Y_1)-(2n_1-2)\eta_1(X_1)\eta_1(Y_1),
						\end{align*} 
						where we have used \eqref{eq:metric} in the third equality. 
						
						For $X_1\in\X(S_1), Y_2\in\X(S_2)$,
						\begin{align*}
							&\Ric^B(X_1,Y_2)=2a[n_1+n_2(a^2+b^2)]\eta_1(X_1)\eta_2(Y_2)-\frac14 \big[ \underbrace{g(2\varphi_1 X_1, 2a \varphi_2 Y_2)}_{=0} \\
							&\quad +\sum_{j=1}^{2n_1} g_1(-2(\eta_1(X_1)\varphi_1 e_j+\Phi_1(e_j,X_1)\xi_1,-2a \eta_2(Y_2)\varphi_1 e_j)\\
							&\quad +\sum_{k=1}^{2n_2} g_2(-2a \eta_1(X_1)\varphi_2 f_k, -2(a^2+b^2)\eta_2(Y_2)\varphi_2 f_k+2\Phi_2(Y_2,f_k)\xi_2)
							\big]\\
							&=2a[n_1+n_2(a^2+b^2)]\eta_1(X_1)\eta_2(Y_2)-\frac14 \big[8 a n_1 \eta_1(X_1)\eta_2(Y_2)+8 n_2 a(a^2+b^2)\eta_1(X_1)\eta_2(Y_2)
							\big ]\\
							&=0.
						\end{align*}
						
						Finally, for $X_2,Y_2\in \X(S_2)$,
						\begin{align*}
							&\Ric^B(X_2,Y_2)=\Ric^2(X_2,Y_2)-2\lambda_{a,b}g_2(X_2,Y_2)+2[n_1 a^2+\lambda_{a,b}+n_2(a^2+b^2)^2-n_2]\eta_2(X_2)\eta_2(Y_2)\\
							&\quad -\frac14 \big[ 4(a^2+b^2) g_2(\varphi_2 X_2, \varphi_2 Y_2)+\sum_{j=1}^{2n_1} g_1(-2a\eta_2(X_2) \varphi_1 e_j, -2a\eta_2(Y_2) \varphi_1 e_j) \\
							&\quad +\sum_{k=1}^{2n_2} \big(g(-2(a^2+b^2)\eta_2(X_2)\varphi_2 f_k,-2(a^2+b^2)\eta_2(Y_2)\varphi_2 f_k) + g(2\Phi_2(X_2,f_k)\xi_2,2\Phi_2(Y_2,f_k)\xi_2)\big)\big]\\
							&=\Ric^2(X_2,Y_2)-2\lambda_{a,b}g_2(X_2,Y_2)+2[n_1 a^2+\lambda_{a,b}+n_2(a^2+b^2)^2-n_2]\eta_2(X_2)\eta_2(Y_2)\\
							&\quad -\frac14 \big[ 4 (a^2+b^2) g_2(\varphi_2 X_2,\varphi_2 Y_2)+8 n_1 a^2 \eta_2(X_2)\eta_2(Y_2)\\
							&\quad +8 n_2(a^2+b^2)^2 \eta_2(X_2)\eta_2(Y_2)+4\sum_{k=1}^{2n_2} (a^2+b^2) g_2(f_k,\varphi_2 X_2)g_2(f_k,\varphi_2 Y_2)\big] \\
							&=\Ric^2(X_2,Y_2)-2\lambda_{a,b}g_2(X_2,Y_2)+2[n_1 a^2+\lambda_{a,b}+n_2(a^2+b^2)^2-n_2]\eta_2(X_2)\eta_2(Y_2)\\
							&\quad-2(a^2+b^2)g_2(\varphi_2 X_2,\varphi_2 Y_2)-2 n_1 a^2 \eta_2(X_2)\eta_2(Y_2)-2n_2(a^2+b^2)^2 \eta_2(X_2)\eta_2(Y_2)\\
							&=\Ric^2(X_2,Y_2)-2(2a^2+2b^2-1)g_2(X_2,Y_2) +[2(2a^2+2b^2-1)-2n_2) ]\eta_2(X_2)\eta_2(Y_2)
						\end{align*}
						where we have used that $g(\xi_2,\xi_2)=a^2+b^2$ in the second equality and \eqref{eq:metric} in the last one.
						
						The last statement concerning $\Ric^B=0$ is clear.
					\end{proof}
					
					\begin{corollary}
						On $M_{0,1}=S_1\times S_2$ equipped with Morimoto's structure $(J_{0,1},g_{0,1})$, the Bismut-Ricci tensor $\Ric^B$ vanishes if and only if $S_1$ and $S_2$ are $\eta$-Einstein with constants 
						\begin{align*}
							(\lambda_1,\nu_1)&=(2,2n_1-2),\\
							(\lambda_2,\nu_2)&=(2,2n_2-2),
						\end{align*} respectively.
					\end{corollary}
					
					\begin{remark}\label{remark:Ric^B(xi)}
						It follows from the expressions for $\Ric^B$ obtained in Theorem \ref{theorem:BRF} that 
						\[ \Ric^B(\xi_1,X) =\Ric^B(\xi_2,X)=0, \] 
						for any $X\in\X(M_{a,b})$.
					\end{remark} 
					
					\smallskip
					
					\begin{remark}
						We note that the name Bismut-Ricci flat has been used recently to denote metric connections with skew-symmetric torsion such that the torsion 3-form is closed and the associated Ricci tensor vanishes (see for instance \cite{GS,PR1,PR2}). In our context, this reduces to Hermitian structures $(J_{a,b},g_{a,b})$ which are SKT (since $dT^B=dd^c \omega=0$) and satisfy $\Ric^B=0$. It follows from Proposition \ref{proposition:skt}  that, in the SKT case, $\dim(S_1\times S_2)\leq 6$ and, in dimension 6, the only Bismut-Ricci flat structures of the form $(J_{a,b},g_{a,b})$ occur on products $S_1\times S_2$ with $\dim S_1=\dim S_2=3$ and $a=0$. Moreover, we will see in \S\ref{section:examples} that  Theorem \ref{theorem:BRF} implies that both $S_1$ and $S_2$ are quotients of the form $\mathbb{S}^3/\Gamma$ with $\Gamma$ a finite subgroup of $\operatorname{SU}(2)$, which is identified with $\mathbb{S}^3$. 
					\end{remark}
					
					\medskip
					
					As an application of Theorem \ref{theorem:BRF}, using a result in \cite{FG}, we obtain information about the canonical bundle of the compact complex manifold $(M_{a,b},J_{a,b})$ when $\Ric^B=0$. 
					
					\begin{proposition}
						The product $M_{a,b}=S_1\times S_2$ of two compact Sasakian manifolds equipped with the Hermitian structure $(J_{a,b},g_{a,b})$ such that $\Ric^B=0$ does not have holomorphically trivial canonical bundle, provided $n_1\geq 1$ and $n_2\geq 1$. 
					\end{proposition}
					
					\begin{proof}
						According to \cite[Theorem 4.1]{FG}, if $(M_{a,b},J_{a,b})$  admitted a non-vanishing holomorphic $(n_1+n_2,0)$-form then the fact that $\Ric^B=0$ would imply that $(M_{a,b},J_{a,b},g_{a,b})$ is conformally balanced. This would mean that $d\theta_{a,b}=0$, however, it follows from \eqref{eq:theta} that \[ d\theta_{a,b}=\frac{4b}{n_1+n_2}(n_2 \Phi_1-n_1 \Phi_2),\] which is non-zero since $n_1\geq 1$ and $n_2\geq 1$. Therefore the canonical bundle of $(M_{a,b},J_{a,b})$ is not trivial.
					\end{proof}
					
					\begin{corollary}
						Let $S$ be a Sasakian $\eta$-Einstein manifold of dimension $\geq 3$ with constants $(\lambda,\nu)$, $\lambda>-2$. Then, $(\mathbb{S}^3\times S,J_{a,b})$ does not have trivial canonical bundle, for $a,b$ such that $a^2+b^2=\frac{\lambda+2}{4}$.
					\end{corollary}
					
					Examples of Sasakian $\eta$-Einstein manifolds with $\lambda>-2$ will appear in \S\ref{section:examples}.
					
					\medskip
					
					We analyze in the next result the CYT condition on $M_{a,b}=S_1\times S_2$, namely, $\rho^B=0$. 
					
					\begin{theorem}\label{theorem:CYT}
						Assuming $n_1\geq1$ and $n_2\geq 1$, the Bismut-Ricci form $\rho^B$ associated to the Hermitian structure $(J_{a,b},g_{a,b})$ on $M_{a,b}=S_1\times S_2$ is given by:
						\begin{align*}
							\rho^B(X_1,Y_1)&=\Ric^1(X_1,\varphi_1 Y_1)-2(2n_1+2an_2-1)\Phi_1(X_1,Y_1),\\
							\rho^B(X_1,Y_2)&=0,\\
							\rho^B(X_2,Y_2)&=\Ric^2(X_2,\varphi_2 Y_2)-2[2an_1+2(a^2+b^2)n_2-1]\Phi_2(X_2,Y_2).
						\end{align*}
						In particular, $(J_{a,b},g_{a,b})$ is CYT if and only if both $S_1$ and $S_2$ are $\eta$-Einstein with constants $(\lambda_1,\nu_1)$ and $(\lambda_2, \nu_2)$ respectively, where
						\begin{equation}\label{eq:lambda-CYT}
							\lambda_1= 4(n_1+an_2)-2,\quad \text{and} \quad 
							\lambda_2= 4(an_1+(a^2+b^2)n_2)-2.
						\end{equation}
					\end{theorem}
					
					\begin{proof} 
						We start by recalling a formula for $\rho^B$ on a $2n$-dimensional Hermitian manifold equipped with the Bismut connection, due to \cite{IP}:
						\[ \rho^B(X,Y)=\Ric^B(X,JY)+(\nabla^B_X \theta) JY+\frac14 \lambda^\omega(X,Y),\]
						where $\theta$ denotes the Lee form and the 2-form $\lambda^\omega$ is defined by \[ \lambda^\omega(X,Y)=\sum_{i=1}^{2n} dT^B(X,Y,u_i,Ju_i),\] in a local orthonormal frame $\{u_i\}_{i=1}^{2n}$. 
						
						In our case of a Sasakian product $M_{a,b}=S_1\times S_2$, we will the use the local orthonormal frame \eqref{eq:frame}. Note that, according to Proposition \ref{proposition:bismut}, we have that $\nabla^B \theta =0$.
						
						We compute first $\lambda^\omega$ (note that these computations only make sense when $n_i\geq 1$, $i=1,2$). From \eqref{eq:T^B} we obtain \[  dT^B=4[\Phi_1^2+2a \Phi_1 \wedge \Phi_2+(a^2+b^2)\Phi_2^2].\]
						This says that $dT^B(X,Y,Z,W)=0$ when one of $X,Y,Z,W$ is $\xi_1$ or $\xi_2$. Thus, 
						\begin{equation}\label{eq:lambdaomega} 
							\iota_{\xi_1} \lambda^\omega=\iota_{\xi_2} 
							\lambda^\omega=0,
						\end{equation} and for $X,Y\in\X(M_{a,b})$,
						\[ \lambda^\omega(X,Y)=\sum_{j=1}^{2n_1} dT^B(X,Y,e_j,\varphi_1 e_j)+\sum_{k=1}^{2n_2} dT^B(X,Y,f_k,\varphi_2 f_k).\]
						To compute these terms we will use a formula for $dT^B $ in terms of $T^B$ that appears in the proof of \cite[Proposition 3.1]{IP}. Given that $\nabla^B T^B=0$, this formula simplifies to
						\begin{align*}
							dT^B(X,Y,Z,W)={\substack{\displaystyle{\mathfrak{S}}\vspace{0,1cm} \\ X,Y,Z}}\; 2 g(T^B(X,Y),T^B(Z,W)), 
						\end{align*}
						where ${\substack{\displaystyle{\mathfrak{S}}\vspace{0,1cm} \\ X,Y,Z}}$ denotes the cyclic sum of $X,Y,Z$.
						Then, it follows from Lemma \ref{lemma:T^B} that, for $X_1,Y_1\in \Gamma(\D_1)$,
						\begin{align*}
							dT^B(X_1,Y_1,e_j,\varphi_1 e_j)&=2\,{\substack{\displaystyle{\mathfrak{S}}\vspace{0,1cm} \\ X_1,Y_1,e_j}}\;g(T^B(X_1,Y_1),T^B(e_j,\varphi_1 e_j))\\
							&=8\,{\substack{\displaystyle{\mathfrak{S}}\vspace{0,1cm} \\ X_1,Y_1,e_j}}\;g(\Phi_1(X_1,Y_1)\xi_1,\Phi_1(e_j,\varphi_1 e_j) \xi_1)\\
							&=8 [\Phi_1(X_1,Y_1)\Phi_1(e_j,\varphi_1 e_j)+\Phi_1(Y_1,e_j)\Phi_1(X_1,\varphi_1 e_j)\\
							&\quad +\Phi_1(e_j,X_1)\Phi_1(Y_1,\varphi_1 e_j)]\\
							&=8[-\Phi_1(X_1,Y_1)+g_1(e_j,\varphi_1 Y_1)g_1(e_j,X_1) -g_1(e_j,\varphi_1 X_1)g_1(e_j,Y_1)]
						\end{align*}
						and
						\begin{align*}
							dT^B(X_1,Y_1,f_k,\varphi_2 f_k)&=2[g(T^B(X_1,Y_1),T^B(f_k,\varphi_2 f_k))+g(T^B(Y_1,f_k),T^B(Y_1,\varphi_2 f_k))\\
							&\quad +g(T^B(f_k,X_1),T^B(Y_1,\varphi_2 f_k))]
						\end{align*}
						It follows from Lemma \ref{lemma:T^B} that $T^B(U_1,U_2)=0$ for $U_i\in \Gamma(\D_i)$, so we arrive at
						\begin{align*}
							dT^B(X_1,Y_1,f_k,\varphi_2 f_k)&=8  g( \Phi_1(X_1,Y_1)\xi_1, \Phi_2(f_k,\varphi_2 f_k)\xi_2) \\
							&= -8a\Phi_1(X_1,Y_1).
						\end{align*}
						Therefore \begin{align*}\lambda^\omega(X_1,Y_1)&=\sum_{j=1}^{2n_1} dT^B(X_1,Y_1,e_j,\varphi_1 e_j)+\sum_{k=1}^{2n_2} dT^B(X_1,Y_1,f_k,\varphi_2 f_k)\\
							&=8[-2n_1 \Phi_1(X_1,Y_1)+2\Phi_1(X_1,Y_1)]-16an_2 \Phi_1(X_1,Y_1)\\
							&=-16(n_1+an_2-1)\Phi_1(X_1,Y_1). \end{align*}
						For $X_1\in \Gamma(\D_1), Y_2\in \Gamma(\D_2)$, it follows from $T^B(U_1,U_2)=0$, when $U_i\in \Gamma(\D_i)$, that \[\lambda^\omega(X_1,Y_2)=0.\]
						Finally, for $X_2,Y_2\in \Gamma(\D_2)$,
						\begin{align*}
							dT^B(X_2,Y_2,e_j,\varphi_1 e_j)&=8 g(\Phi_2(X_2,Y_2)\xi_2,\Phi_1(e_j,\varphi_1 e_j)\xi_1)\\
							&=8a\Phi_2(X_2,Y_2)\Phi_1(e_j,\varphi_1 e_j)\\
							&=-8a  \Phi_2(X_2,Y_2),
						\end{align*}
						and
						\begin{align*}
							dT^B(X_2,Y_2,f_k,\varphi_2 f_k)&=2\,{\substack{\displaystyle{\mathfrak{S}}\vspace{0,1cm} \\ X_2,Y_2,f_k}}\; g(\Phi_2(X_2,Y_2)\xi_2,\Phi_2(f_k,\varphi_2 f_k)\xi_2)\\
							&=8(a^2+b^2) [\Phi_2(X_2,Y_2)\Phi_2(f_k,\varphi_2 f_k)+\Phi_2(Y_2,f_k)\Phi_2(X_2,\varphi_2 f_k)\\
							&\quad+\Phi_2(f_k,X_2)\Phi_2(Y_2,\varphi_2 f_k)]\\
							&=8(a^2+b^2)[-\Phi_2(X_2,Y_2)+g_2(f_k,\varphi_2 Y_2)g_2(f_k,X_2)\\
							&\quad-g_2(f_k,\varphi_2 X_2)g_2(f_k,Y_2)].
						\end{align*}
						Hence, 
						\begin{align*}
							\lambda^\omega(X_2,Y_2)&=\sum_{j=1}^{2n_1} dT^B(X_2,Y_2,e_j,\varphi_1 e_j)+\sum_{k=1}^{2n_2} dT^B(X_2,Y_2,f_k,\varphi_2 f_k)\\
							&=-16an_1\Phi_2(X_2,Y_2)+8(a^2+b^2)[-2n_2\Phi_2(X_2,Y_2)+2\Phi_2(X_2,Y_2)]\\
							&=-16[an_1+(a^2+b^2)(n_2-1)]\Phi_2(X_2,Y_2).
						\end{align*}
						Now we proceed to compute $\rho^B$. First, note that $\rho^B(X_1,Y_2)=0$ when $X_1\in\X(S_1)$ and $Y_2\in\X(S_2)$ since $\Ric^B(X_1,Y_2)=\lambda^\omega(X_1,Y_2)=0$. Moreover, it follows easily from Remark \ref{remark:Ric^B(xi)} and \eqref{eq:lambdaomega} that \begin{equation}\label{eq:rho-xi}
							\rho^B(\xi_i,X)=0, \quad X\in\X(M_{a,b}).
						\end{equation}  Therefore, it is enough to compute $\rho^B(X_i,Y_i)$, $i=1,2$, for $X_i,Y_i\in\Gamma(\D_i)$.
						
						From the computations above we obtain that, for $X_1,Y_1\in \Gamma(\D_1)$,
						\begin{align*}
							\rho^B(X_1,Y_1)&=\Ric^B(X_1,\varphi_1 Y_1)+\frac{1}{4} \lambda^\omega(X_1,Y_1)\\
							&=\Ric^1(X_1,\varphi_1 Y_1)-2\Phi_1(X_1,Y_1)-4(n_1+an_2-1)\Phi_1(X_1,Y_1)\\
							&=\Ric^1(X_1,\varphi_1 Y_1)-2(2n_1+2an_2-1)\Phi_1(X_1,Y_1),
						\end{align*}
						where we have used Theorem \ref{theorem:BRF} in the second equality.
						
						For $X_2,Y_2\in\Gamma(\D_2)$,
						\begin{align*}
							\rho^B(X_2,Y_2)&=\Ric^B(X_2,\varphi_2 Y_2)-4(an_1+(a^2+b^2)(n_2-1))\Phi_2(X_2,Y_2)\\
							&=\Ric^2(X_2,\varphi_2 Y_2)-2(2a^2+2b^2-1)\Phi_2(X_2,Y_2)\\
							&\quad -4(an_1+(a^2+b^2)(n_2-1))\Phi_2(X_2,Y_2)\\
							&=\Ric^2(X_2,\varphi_2 Y_2)-2[2an_1+2(a^2+b^2)n_2-1]\Phi_2(X_2,Y_2),
						\end{align*}
						where we have used again Theorem \ref{theorem:BRF} in the second equality.
						
						Therefore, according to \eqref{eq:Ric-D} and using that $\varphi_i$ is an isomorphism on $\mathcal{D}_i$, $i=1,2$, we arrive at
						\[ \rho^B\equiv 0 \iff \begin{cases}
							\Ric^1=\lambda_1 g_1 + (2n_1-\lambda_1) \eta_1 \otimes \eta_1,\\
							\Ric^2=\lambda_2 g_2 + (2n_2-\lambda_2) \eta_2 \otimes \eta_2,
						\end{cases}\] where $\lambda_1=4(n_1+an_2)-2$ and $\lambda_2=4(an_1+(a^2+b^2)n_2)-2$, and this finishes the proof.
					\end{proof}
					
					\begin{corollary} 
						Assuming $n_1\geq1, n_2\geq 1$, the product manifold $M_{0,1}=S_1\times S_2$ equipped with Morimoto's structure $(J_{0,1},g_{0,1})$ is CYT if and only if $S_1$ and $S_2$ are $\eta$-Einstein with constants $(\lambda_1,\nu_1)$ and $(\lambda_2,\nu_2)$ respectively, where
						\[\lambda_1=4n_1-2,\quad \text{and} \quad 
						\lambda_2=4n_2-2\]
					\end{corollary}
					
					\begin{remark}\label{remark:missing cases}
						We analyze here the missing cases $n_1=0$ or $n_2=0$ (with $n_1+n_2\geq 1$). 
						Following the lines of the proof of Theorem \ref{theorem:CYT} we get:
						\begin{itemize}
							\item When $n_2=0$, the Hermitian structure $(J_{a,b},g_{a,b})$ is CYT if and only if $S_1$ is $\eta$-Einstein with constants $(\lambda_1,\nu_1)$, $\lambda_1=4n_1-2$.
							\item When $n_1=0$, the Hermitian structure $(J_{a,b},g_{a,b})$ is CYT if and only if $S_2$ is $\eta$-Einstein with constants $(\lambda_2,\nu_2)$, $\lambda_2=4(a^2+b^2)n_2-2$.
						\end{itemize}
					\end{remark}
					
					\medskip
					
					Using the expression for the Bismut-Ricci form $\rho^B$ obtained in Theorem \ref{theorem:CYT}. we can determine when the Hermitian structure $(J_{a,b},g_{a,b})$ on a Sasakian product $M_{a,b}=S_1\times S_2$ is \textit{static}. This notion was introduced by Streets and Tian in \cite{ST1,ST2}: an SKT Hermitian metric $g$ on a complex manifold $(M^{2n},J)$ is called \textit{static} if its Bismut-Ricci form satisfies
					\begin{equation}\label{eq:static}
						(\rho^B)^{1,1}=\alpha \omega,\quad \alpha \in \R,
					\end{equation}
					where $(\rho^B)^{1,1}$ denotes the $(1,1)$-component of $\rho^B$ given by $(\rho^B)^{1,1}(\cdot,\cdot)=\frac12 (\rho^B(\cdot,\cdot)+\rho^B(J\cdot,J\cdot))$. Static metrics are closely related to the \textit{pluriclosed flow}, introduced in \cite{ST1}, which is the parabolic flow for SKT metrics defined by
					\[ \frac{\partial}{\partial t} \omega = -(\rho^B)^{1,1}, \quad \omega(0)=\omega_0
					.\]
					Thus, static metrics are to the pluriclosed flow what Einstein metrics are to the Ricci flow. Furthermore, when $\alpha=0$ these metrics are fixed points of the pluriclosed flow.
					
					The following relation between $\rho^B(JX,JY)$ and $\rho^B(X,Y)$ was proved  in \cite[Corollary 3.2]{IP}:
					\begin{equation}\label{eq:rho^B=(rho^B){1,1}}
						\rho^B(JX,JY)-\rho(X,Y)=\delta T^B(JX,Y)-(\nabla^B_{JX} \theta)Y+(\nabla^B_Y \theta)JX.
					\end{equation}
					Using \eqref{eq:codif T^B}, \eqref{eq:rho^B=(rho^B){1,1}} and $\nabla^B \theta=0$ (see Proposition \ref{proposition:bismut}) we obtain that $\rho^B(JX,JY)=\rho^B(X,Y)$ and as a consequence,
					\begin{equation}\label{eq:rhoiguales}
						(\rho^B)^{1,1}=\rho^B,
					\end{equation}
					that is, $\rho^B$ is $J$-invariant. Due to \eqref{eq:rhoiguales} and taking \eqref{eq:omega} into account, condition \eqref{eq:static} becomes \[ \rho^B=\alpha \,\omega_{a,b}=\alpha(\Phi_1+\Phi_2-b\eta_1\wedge \eta_2).\]
					Since $\rho^B(\xi_1,\xi_2)=0$ (due to \eqref{eq:rho-xi}) and $\omega_{a,b}(\xi_1,\xi_2)=-b\neq 0$, we obtain that $\rho^B(\xi_1,\xi_2)=\alpha \,\omega_{a,b}(\xi_1,\xi_2)$ if and only if $\alpha=0$. Therefore, condition \eqref{eq:static} reduces to the CYT condition. To sum up, the Hermitian structure $(J_{a,b},g_{a,b})$ on a product of Sasakian manifolds $M_{a,b}=S_1\times S_2$ satisfies \eqref{eq:static} if and only if $\alpha=0$ and $(J_{a,b},g_{a,b})$ is CYT. 
					
					Let us recall that $(J_{a,b},g_{a,b})$ is SKT only in dimensions 4 and 6: any such structure is SKT in dimension $4$, and in dimension 6, $\dim S_1=\dim S_2=3$ and $a=0$ (see Proposition \ref{proposition:skt}). In dimension 4, we will see 
					in \S\ref{section:examples} that Remark \ref{remark:missing cases} implies that one factor is one-dimensional and the other one is a quotient of the form $\mathbb{S}^3/\Gamma$ with $\Gamma$ a finite subgroup of $\operatorname{SU}(2)\simeq \mathbb{S}^3$. In dimension 6, we will see in \S\ref{section:examples} that Theorem \ref{theorem:CYT} implies that both $S_1$ and $S_2$ are quotients of this form. To sum up,
					
					\begin{proposition}
						The Hermitian structure $(J_{a,b},g_{a,b})$ on $M_{a,b}=S_1\times S_2$ is static if and only if
						\begin{enumerate}
							\item [$(a)$] $\dim M_{a,b}=4$, one of the Sasakian factors is 1-dimensional and the other one is a quotient of $\mathbb{S}^3$ by a finite subgroup, or
							\item [$(b)$] $\dim S_1=\dim S_2=3$, $a=0$ and both $S_1$ and $S_2$ are quotients of $\mathbb{S}^3$ by finite subgroups.
						\end{enumerate}
						Moreover, these metrics are fixed points of the pluriclosed flow.
					\end{proposition}
					
					\begin{remark}
						The fact that $\mathbb{S}^3\times \mathbb{S}^3$ carries a CYT structure is well known (see \cite{GGP}).
					\end{remark}
					
					\medskip 
					
					\subsection{Examples}\label{section:examples}
					
					In order to exhibit examples of Sasakian $\eta$-Einstein manifolds that fulfill the conditions in Theorems \ref{theorem:BRF} and \ref{theorem:CYT}, we recall the notion of $\mathcal{D}$-homothetic deformations (or simply $\mathcal{D}$-homotheties), introduced by Tanno in \cite{Tan}. Given a Sasakian manifold $(S,\varphi,\xi,\eta,g)$, consider the transformation 
					\[ \varphi'=\varphi, \quad \xi'=s^{-1}\xi, \quad \eta'=s\eta, \quad  g'=sg+s(s-1)\eta\otimes \eta, \]
					for any real constant $s>0$. Then $(\varphi',\xi',\eta',g')$ is again a Sasakian structure on $S$. Moreover, in the case of Sasakian $\eta$-Einstein manifolds, there is the following result:
					
					\begin{proposition}\cite[Proposition 18]{BGM}\label{proposition:deformation}
						Let $(S,\varphi,\xi,\eta,g)$ be a $(2n+1)$-dimensional Sasakian $\eta$-Einstein manifold with constants $(\lambda,\nu)$, and consider a $\D$-homothetic structure $(\varphi',\xi',\eta',g')$ as above. Then, $(S,\varphi',\xi',\eta',g')$ is also $\eta$-Einstein with constants \[ \lambda'=\frac{\lambda+2-2s}{s},\qquad \nu'=2n-\frac{\lambda+2-2s}{s}.\]
					\end{proposition}
					
					In this article we will use the following terminology: a Sasakian $\eta$-Einstein manifold $S$  will be called\footnote{The notions of positive, negative and null are defined in \cite{BGM} for general Sasakian manifolds in terms of the basic first Chern class, and they reduce to the stated inequalities for $\lambda$ in the case of $\eta$-Einstein manifolds.} \textit{positive} if $\lambda>-2$, \textit{null} if $\lambda=-2$ and \textit{negative} if $\lambda<-2$. It follows from Proposition \ref{proposition:deformation} that $\D$-homotheties preserve positive, null and negative $\eta$-Einstein manifolds, respectively. Positive $\eta$-Einstein manifolds include the well-known family of Sasaki-Einstein manifolds\footnote{The literature on Sasaki-Einstein metrics is vast, for instance the whole Chapter 11 of \cite{BG} is devoted to Sasaki-Einstein metrics (see also \cite{Spa}).}, which is precisely the case when $\lambda=\dim S-1$ and $\nu=0$. For instance, the odd-dimensional spheres with the standard Sasakian structure are Einstein (in fact, it was recently proved in \cite{LST} that there are infinitely many families of Sasaki-Einstein metrics on every odd-dimensional standard sphere of dimension at least 5). It follows from the Bonnet-Myers theorem that a manifold admitting a positive $\eta$-Einstein structure is compact and has finite fundamental group.
					
					\medskip 
					
					We can rephrase Theorems \ref{theorem:BRF} and  \ref{theorem:CYT} in terms of $S_1$ and $S_2$ being positive, null or negative Sasakian $\eta$-Einstein.
					Since $\lambda_1=2$ and $\lambda_2=2(2a^2+2b^2)-2>-2$, we obtain
					
					\begin{theorem}\label{theorem:BRFreph}
						Let $M=S_1\times S_2$ be the product of two Sasakian manifolds. Then, after possibly applying a $\D$-homothety to each Sasakian structure, $M$ admits a Hermitian structure of the form $(J_{a,b},g_{a,b})$ for some $a,b\in\R$, $b\neq 0$, such that $\Ric^B=0$ if and only if $S_1$ and $S_2$ admit Sasaki-Einstein metrics. 
					\end{theorem}
					
					To rephrase Theorem \ref{theorem:CYT} let us analyze when \[ \lambda_1=4(n_1+an_2)-2 \quad  \text{and} \quad  \lambda_2=4(an_1+(a^2+b^2)n_2)-2 \] are $<-2, =-2$ or $>-2$, respectively. Note that 
					\[ \begin{cases}
						\lambda_1 \geq -2 \iff n_1+an_2\geq 0 \quad \text{and}\\
						\lambda_2 \geq -2 \iff an_1+(a^2+b^2)n_2 \geq 0,
					\end{cases} 
					\]
					Assuming $n_1\geq 1$ and $n_2\geq 1$ we analyze all the combinations for $\lambda_1$ and $\lambda_2$.
					
					Case (i): When $\lambda_1=-2$ we obtain that $a=-\frac{n_1}{n_2}$ and thus $\lambda_2=4 b^2 n_2-2$. Given that $b\neq 0$, it is clear that $\lambda_2>-2$, and in this case any $b\neq 0$ works.
					
					Case (ii):  When $\lambda_1>-2$, $\lambda_2$ can be $<-2$, $=-2$ or $>-2$. Indeed, by possibly performing a $\D$-homothety to $S_1$ we may assume that $\lambda_1=2n_1-2$, so $a=-\frac{n_1}{2n_2}$ and thus $\lambda_2=\frac{4b^2 n_2^2-n_1^2}{n_2}-2$. If $\lambda_2=-2$, we can choose  $b=\frac{n_1}{2n_2}$. If $\lambda_2>-2$, after possibly performing a $\D$-homothety to $S_2$ we may assume that $\lambda_2=\frac{3n_1^2}{n_2}-2$ and choose $b=\frac{n_1}{n_2}$. Analogously, if $\lambda_2<-2$, we may assume that $\lambda_2=-\frac{3n_1^2}{4n_2}-2$ and choose $b=\frac{n_1}{4n_2}$.
					
					
					Case (iii): When $\lambda_1<-2$, we have that 
					$a<-\frac{n_1}{n_2}$. Then we have to determine the sign of \[P(a):=an_1+(a^2+b^2)n_2=n_2\left[\left(a+\frac{n_1}{2n_2}\right)^2+b^2-\frac{n_1^2}{4n_2^2}\right].\] If $b^2-\frac{n_1^2}{4n_2^2}\geq 0$ then $P(a)>0$ provided $a\neq -\frac{n_1}{2n_2}$. In particular $P(a)>0$ for $a<-\frac{n_1}{n_2}$.
					
					If $b^2-\frac{n_1^2}{4n_2^2}<0$, then $P(a)\leq 0$ if and only if $\frac{-n_1-\sqrt{n_1^2-4b^2n_2^2}}{2n_2}\leq a\leq \frac{-n_1+\sqrt{n_1^2-4b^2n_2^2}}{2n_2}$. Since $-\frac{n_1}{n_2}<\frac{-n_1-\sqrt{n_1^2-4b^2n_2^2}}{2n_2}$, we arrive at $P(a)>0$ for $a<-\frac{n_1}{n_2}$. To sum up, $\lambda_1<-2$ implies $\lambda_2>-2$. In this case, after possibly applying a $\mathcal{D}$-homothety to $S_1$ and $S_2$ we can assume that $\lambda_1=-n_1-2$ and $\lambda_2=\frac{3 n_1^2}{2 n_2}-2$. Then, it is easily seen that $a=-\frac{5n_1}{4n_2}$ and $b=\frac{n_1}{4n_2}$ work.
					
					\medskip 
					
					Now we can rephrase Theorem \ref{theorem:CYT} as follows.
					
					\begin{theorem}\label{theorem:CYTreph}
						Let $M=S_1\times S_2$ be the product of two Sasakian manifolds and assume $n_1\geq1$, $n_2\geq 1$. Then, after possibly applying a $\D$-homothety to each Sasakian structure, $M$ admits a CYT Hermitian structure of the form $(J_{a,b},g_{a,b})$ for some $a,b\in\R$, $b\neq 0$ if and only if $S_1$ and $S_2$ are $\eta$-Einstein and one of the following holds:
						\begin{itemize} 
							\item[$\ri$] $S_1$ is positive and $S_2$ is arbitrary, 
							\item[$\rii$] $S_1$ is negative or null, and $S_2$ is positive.
						\end{itemize}
					\end{theorem}
					
					\begin{remark}
						According to Remark \ref{remark:missing cases}, when one of the Sasakian factors is one-dimensional, the CYT condition reduces to the other $\eta$-Einstein factor being positive. 
					\end{remark}
					
					\begin{remark} 
						Since odd-dimensional spheres equipped with their usual structure are Sasaki-Einstein, Theorem \ref{theorem:CYTreph} shows the existence of CYT Hermitian structures on Calabi-Eckmann manifolds. Note that this was already proved in \cite[Corollary 4.8]{Bar} where, more generally, the existence of CYT structures on principal bundles over Hermitian manifolds with complex tori as fibers is proved.
					\end{remark}
					
					\medskip
					
					We point out that if we begin with two Sasaki-Einstein manifolds there is no need to perform any $\mathcal D$-homothety on the factors in order to obtain a CYT structure on their product. Indeed, we have the following result:
					
					\begin{proposition}
						Let $S_1$ and $S_2$ be two Sasaki-Einstein manifolds with $\dim S_i = 2n_i+1$, $n_i\geq 1$, for $i=1,2$. Then $S_1\times S_2$ admits a CYT structure $(J_{a,b},g_{a,b})$ with 
						\[ a=-\frac{n_1-1}{2n_2}, \quad b^2= \frac{(n_1-1)(n_1+1)+2n_2(n_2+1)}{4n_2^2}.\]
					\end{proposition}
					
					\begin{proof}
						The proof follows by solving for $a$ and $b$ in \eqref{eq:lambda-CYT}, with $\lambda_1=2n_1$ and $\lambda_2=2n_2$.
					\end{proof}
					
					\medskip 
					
					\begin{example}[Left invariant $\eta$-Einstein structures on Lie groups]
						
						According to \cite{Geiges}, a 3-dimensional compact Sasakian manifold is diffeomorphic to $\mathbb{S}^3/\Gamma$, $H_3/\Gamma$ or $\widetilde{\operatorname{SL}}(2,\R)/\Gamma$, where in each case $\Gamma$ is a uniform lattice (i.e., a co-compact discrete subgroup). 
						It is known that these 3 model geometries correspond precisely to positive, null or negative Sasakian $\eta$-Einstein structures. In the table below we show an explicit $\eta$-Einstein structure on the corresponding Lie algebras $\sl(2,\R)$, $\h_3$ and $\su(2)$, all of them spanned by an orthonormal  basis $\{e_1,e_2,e_3\}$.
						
						\begin{table}[H]
							\centering
							\begin{tabular}{|c|c|c|c|c|}
								\hline
								Lie algebra & Lie brackets & Sasakian structure & $\lambda$ \\
								\hline
								$\su(2)$ & $[e_1,e_2]=2e_3, [e_2,e_3]=2e_1, [e_3,e_1]=2e_2$ & $\xi=e_3, \eta=e^3, \varphi e_1=e_2$ &2 \\
								$\h_3$ & $[e_1,e_2]=2e_3$ & $ \xi=e_3, \eta=e^3, \varphi e_1=e_2$ & $-2$  \\
								$\sl(2,\R)$ & $[e_1,e_2]=2e_3, [e_2,e_3]=-e_1, [e_3,e_1]=-e_2$ & $\xi=e_3$, $\eta=e^3$, $\varphi e_1=e_2$ & $-4$ \\
								\hline
							\end{tabular}
							\caption{3-dimensional Sasakian Lie algebras}
							\label{table:Sasakian 3 Lie algebras}
						\end{table}
						Using Theorems \ref{theorem:BRF} and \ref{theorem:CYT} we recover the Hermitian structure $(J_{0,1},g_{0,1})$ on $\su(2)\times\su(2)$, which has the special feature that $\nabla^B\equiv 0$ on left invariant vector fields, so it is Bismut flat and in particular $\Ric^B=0$ and $\rho^B=0$. Moreover, we obtain the following CYT Hermitian structures $(J_{a,b},g_{a,b})$:
						\begin{itemize}
							\item on $\h_3\times \su(2)$, with $(a,b)=(-1, 1)$,
							\item on $\sl(2,\R)\times \su(2)$, with $(a,b)=(-\frac{3}{2}, \frac{1}{2})$,
							\item on $\su(2)\times \h_3$ (after applying a $\D$-homothety to $\su(2)$ so that $\lambda_1=0$), with $(a,b)=(-\frac12, \frac12)$,
							\item on $\su(2)\times \sl(2,\R)$ (after applying a $\D$-homothety to $\su(2)$ and $\sl(2,\R)$ so that $\lambda_1=0$, $\lambda_2=-\frac{11}{4}$), with $(a,b)=(-\frac{1}{2},\frac{1}{4})$.
						\end{itemize} 
						
						\medskip
						
						The Sasakian manifolds obtained as quotients of  $\operatorname{SU}(2)$, $\widetilde{\operatorname{SL}}(2,\R)$ and $H_3$ by a uniform lattice carry induced $\eta$-Einstein structures with the same constant $\lambda$, so that their products admit induced structures such that $\Ric^B=0$ (in the case of $\operatorname{SU}(2)\times \operatorname{SU}(2)$) and CYT structures (in all the cases above).
						
						In higher dimensions, a family of null $\eta$-Einstein Lie algebras is given by the $(2n+1)$-dimensional Heisenberg Lie algebra $\h_{2n+1}$. In fact, it can be seen in the same way as for $\h_3$ that $\h_{2n+1}$, spanned by $\{X_1,\ldots,X_{2n},\xi\}$ with brackets $[X_{2i-1},X_{2i}]=2\xi$ for all $1\leq i\leq n$, is an example of a null Sasakian Lie algebra with $\lambda=-2$. Therefore $\h_{2n+1}\times \su(2)$ admits a CYT structure given by $(J_{a,b},g_{a,b})$ with $(a,b)=(-n,1)$, while $\su(2)\times \h_{2n+1}$ admits a CYT structure given by $(a,b)=(-\frac{1}{2n}, \frac{1}{2n})$ (after applying a $\D$-homothety to $\su(2)$ so that $\lambda_1=0$). The associated simply-connected Lie groups $H_{2n+1}\times \operatorname{SU}(2)$ and $\operatorname{SU}(2)\times H_{2n+1}$ admit uniform lattices, so that we obtain compact CYT manifolds.
						
						In contrast, there is just one example of a positive $\eta$-Einstein Lie algebra, which is given by $\su(2)$. Indeed, a Lie group $G$ admitting a left invariant positive $\eta$-Einstein  structure (or equivalently, a Sasaki-Einstein structure) is compact and has finite fundamental group, therefore $G$ is semisimple. According to \cite{BW}, the only semisimple Lie algebras carrying a contact form are $\mathfrak{su}(2)$ and $\mathfrak{sl}(2,\R)$; hence the Lie algebra of $G$ is $\mathfrak{su}(2)$. 
						
						We do not know yet of other examples, besides $\widetilde{\operatorname{SL}}(2,\R)$, of left invariant negative $\eta$-Einstein structures on Lie groups which admit lattices.
						
						Note that by a result in \cite{AFV}, the only 5-dimensional Sasakian $\eta$-Einstein Lie algebra whose associated simply-connected Lie group admits lattices is $\h_5$ (and this is null).
					\end{example}
					
					\medskip

					\begin{example}[Positive $\eta$-Einstein structures]
						Tanno was the first to observe that by applying a suitable $\D$-homothety in the positive case one obtains a Sasaki-Einstein structure, and he used this to prove that the unit tangent bundle of $\mathbb{S}^n$ has a homogeneous Sasaki-Einstein structure \cite{Tan2}. In particular, $\mathbb{S}^2\times \mathbb{S}^3$ has a Sasaki-Einstein structure. Over the last decades infinitely many  Sasaki-Einstein structures were shown to exist 
						on connected sums of $\mathbb{S}^2\times \mathbb{S}^3$ and also on 5-dimensional manifolds which are not connected sums of $\mathbb{S}^2 \times \mathbb{S}^3$, including infinitely many rational homology 5-spheres. Similar results hold also in higher dimensions (see \cite{Spa} and references therein).  
						
						A special class of Sasaki-Einstein manifolds is given by \textit{3-Sasakian manifolds} which are those Riemannian manifolds whose metric cone is hyperK\"ahler. In particular, a $3$-Sasakian manifold has dimension $4n+3$, with $n\geq 0$, and it carries $3$ different compatible Sasakian structures. They were introduced by C. Udri\c{s}te \cite{Ud} and Y. Kuo \cite{Kuo} in 1969 and 1970, respectively. According to Theorem \ref{theorem:BRFreph} the product of any two of these manifolds can be endowed with a Hermitian structure $(J_{a,b},g_{a,b})$ satisfying $\Ric^B=0$ or $\rho^B=0$, choosing properly the values of $a$ and $b\neq 0$, after possibly applying a $\D$-homothety.
					\end{example}
					
					\smallskip
					
					To finish this section we briefly mention another construction which furnishes many examples of $\eta$-Einstein metrics, especially negative and null. According to \cite{BGM}, many interesting Sasakian examples can be found on links of isolated hypersurface singularities. Moreover, some of them carry $\eta$-Einstein structures. For more details on this construction see for instance \cite[Chapter 9]{BG}, \cite[Section 6]{BGM} and \cite[Section 3.4]{Spa}.
					
					\medskip
					
					\begin{example} [$\eta$-Einstein structures on links]
						Consider the affine space $\C^{n+1}$ together with a weighted $\C^\star$-action given by \[ (z_0,\ldots,z_n)\mapsto (\lambda^{w_0} z_0, \ldots, \lambda^{w_n} z_n),\] where the weights $w_j$ are positive integers such that $\operatorname{gcd}(w_0,\ldots,w_n)=1$. A weighted homogeneous polynomial with weights $w=(w_0,\ldots,w_n)\in \N^{n+1}$ of degree $d$ is a polynomial $f\in \C[z_0,\ldots,z_n]$ such that \[ f(\lambda^{w_0} z_0, \ldots, \lambda^{w_n} z_n)=\lambda^d f(z_0,\ldots,z_n).\]
						Assume that the origin is an isolated singularity of $\{f=0\}$. Then, the \textit{link} of $f$ is defined by \[ L_f=\{f=0\} \cap \mathbb{S}^{2n+1},\] where $\mathbb{S}^{2n+1}$ is the unit sphere in $\C^{n+1}$, and it is a smooth manifold of dimension $2n-1$ which by the Milnor Fibration Theorem is $(n-2)$-connected. The link $L_f$ is endowed with a natural Sasakian structure $S_{w,f}=S_w|_{L_f}$ inherited as a Sasakian submanifold of $\mathbb{S}^{2n+1}$ with its weighted Sasakian structure $(\Phi_w,\xi_w,\eta_w,g_w)$ which in the standard coordinates on $\C^{n+1}\equiv \R^{2n+2}$ is determined by
						\[  \eta_w=\frac{\sum_{i=0}^n (x_i dy_i-y_i dx_i)}{\sum_{i=0}^n w_i(x_i^2+y_i^2)}, \quad \xi_w=\sum_{i=0}^n w_i (x_i \del y_i-y_i \del x_i).\] 
						
						Regarding the existence of $\eta$-Einstein on these links, there is the following result in \cite{BGK} that establishes the existence of negative or null $\eta$-Einstein structures.
						
						\begin{theorem}\cite{BGK}\label{theorem:BGK}
							Let $f$ be a non-degenerate weighted homogeneous polynomial of degree $d$ and weight vector $w$, and let $|w| = w_0 + \cdots + w_n$. Consider the induced Sasakian structure $S_{w,f}$ on the link $L_f$. 
							
							\begin{enumerate} 
								\item If $|w|=d$, then there exists a null $\eta$-Einstein structure on $L_f$,
								\item If $|w|<d$, then there exists a negative $\eta$-Einstein structure on $L_f$.
							\end{enumerate}
							The $\eta$-Einstein structures are obtained by deforming suitably the induced Sasakian structure $S_{w,f}$.
						\end{theorem}
						
						\begin{remark}
							In the case $|w|>d$, there are obstructions for the existence of positive $\eta$-Einstein structures on $L_f$.
						\end{remark}
						
						One well-known example is the \textit{Brieskorn-Pham} link $L(a_0,\ldots,a_n)$ associated to the polynomial $f(z)=z_0^{a_0}+\cdots+z_n^{a_n}$. It can be seen that the weighted degree of $f$ is $d=\operatorname{lcm}(a_0,\ldots,a_n)$ and the weights are $w_j=\frac{d}{a_j}$. Then by Theorem \ref{theorem:BGK}, 
						\begin{itemize} 
							\item when $\sum_i \frac{1}{a_i}=1$ there is a null $\eta$-Einstein structure on $L(a_0,\ldots,a_n)$,
							\item when $\sum_i \frac{1}{a_i}<1$ there is a negative $\eta$-Einstein structure on $L(a_0,\ldots,a_n)$.
						\end{itemize}
						
						According to Theorem \ref{theorem:CYTreph}, the product $\mathbb{S}^{2m+1}\times L(a_0,\ldots,a_n)$ admits CYT structures of the form $(J_{a,b},g_{a,b})$ as long as $\sum_{i} \frac{1}{a_i}\leq 1$, by possibly applying a $\D$-homothety.
					\end{example}
					
					\medskip
					
					\section{Appendix}
					
					As mentioned before, we provide here a proof of Proposition \ref{proposition:integrable-wood}, for the sake of completeness.
					
					\begin{proof}[Proof of Proposition \ref{proposition:integrable-wood}]
						Let $\{e_1,\ldots,e_{2n}\}$ be an orthonormal local frame satisfying $Je_{2i-1}=e_{2i}$ for $1\leq i\leq n$. Using this frame we compute by definition both sides of the equality we want to prove. For $X\in\X(M)$, using that $J(\nabla_U J)=-(\nabla_U J)J$ for all $U\in\X(M)$, we have 
						\begin{align*}
							[J,\nabla^*\nabla J](X)&=\sum_{i=1}^{2n} \underbrace{J(\nabla_{e_i} (\nabla_{e_i} J))(X)}_{\Circled{1}}-\underbrace{(\nabla_{e_i} (\nabla_{e_i} J))(JX)}_{\Circled{2}}-2J(\nabla_{\nabla_{e_i} e_i} J)(X),
						\end{align*}
						Recall that the integrability of $J$ is equivalent to $\nabla_{JU} J=J(\nabla_U J)$ for any $U\in\X(M)$. Using this fact when writing $(\nabla_{e_i} J)=-(\nabla_{J^2 e_i} J)$, we obtain
						\begin{align*} \Circled{1}&=-J(\nabla_{e_i} J(\nabla_{Je_i} J))(X)\\
							&=-J\nabla_{e_i} (J(\nabla_{Je_i} J)X)-(\nabla_{Je_i} J)(\nabla_{e_i} X)\\
							&=-J\nabla_{e_i} J\nabla_{Je_i} JX-J\nabla_{e_i} \nabla_{Je_i} X-\nabla_{Je_i} J\nabla_{e_i} X+J\nabla_{Je_i} \nabla_{e_i} X\\
							&=-J\nabla_{e_i} J\nabla_{Je_i} JX-\nabla_{Je_i} J\nabla_{e_i} X-J R(e_i,Je_i)X-J\nabla_{[e_i,Je_i]} X,
						\end{align*}
						\begin{align*}
							\Circled{2}&=-(\nabla_{e_i} J(\nabla_{Je_i} J))(JX)\\
							&=-\nabla_{e_i} (J(\nabla_{Je_i} J)(JX))+J(\nabla_{Je_i} J)(\nabla_{e_i} JX)\\
							&=\nabla_{e_i} J \nabla_{Je_i} X-\nabla_{e_i} \nabla_{Je_i} JX+J\nabla_{Je_i} J\nabla_{e_i} JX+\nabla_{Je_i} \nabla_{e_i} JX\\
							&=\nabla_{e_i} J \nabla_{Je_i} X+J\nabla_{Je_i} J\nabla_{e_i} JX-R(e_i,Je_i)JX-\nabla_{[e_i,Je_i]} JX. 
						\end{align*}
						
						Hence, 
						\begin{align*} 
							[J,\nabla^* \nabla J](X)&=-2[J,P](X)+\sum_{i=1}^{2n} (\nabla_{[e_i,Je_i]} J)X-2\sum_{i=1}^{2n} J(\nabla_{\nabla_{e_i} e_i} J)(X)\\
							&\quad-\sum_{i=1}^{2n} (J\nabla_{e_i} J\nabla_{Je_i} JX+\nabla_{Je_i} J\nabla_{e_i} X+\nabla_{e_i} J \nabla_{Je_i} X+J\nabla_{Je_i} J\nabla_{e_i} JX).\end{align*}
						Note that in the chosen  $J$-adapted frame $\{e_i\}$ the last sum equals zero, since replacing $e_i$ by $Je_i$ gives the same terms with opposite sign. Thus, 
						\begin{equation}\label{eq:wood-1}
							[J,\nabla^* \nabla J](X) =-2[J,P](X)+\sum_{i=1}^{2n} (\nabla_{[e_i,Je_i]} J)X-2\sum_{i=1}^{2n} J(\nabla_{\nabla_{e_i} e_i} J)(X). 
						\end{equation} 
						Now, using \eqref{eq:delta-J} with our $J$-adapted frame we get
						\begin{align*} 
							2(\nabla_{\delta J} J)(X)&=2\sum_{i=1}^{2n} (\nabla_{(\nabla_{e_i} J)e_i} J)(X)\\
							&=2\sum_{i=1}^{2n} (\nabla_{\nabla_{e_i} Je_i} J)(X)-2\sum_{i=1}^{2n} (\nabla_{J \nabla_{e_i} e_i} J)(X)\\
							&=\sum_{i=1}^{2n} (\nabla_{\nabla_{e_i} Je_i} J+\nabla_{\nabla_{Je_i} e_i} J)(X)+\sum_{i=1}^{2n} (\nabla_{[e_i, Je_i]} J)(X) -2\sum_{i=1}^{2n} J(\nabla_{ \nabla_{e_i} e_i} J)(X).
						\end{align*}
						Replacing $e_i$ by $Je_i$ in the first sum we obtain the same terms with opposite sign, and thus this sum equals zero. Therefore,
						\begin{equation}\label{eq:wood-2}
							2(\nabla_{\delta J} J)(X) = \sum_{i=1}^{2n} (\nabla_{[e_i, Je_i]} J)(X)-2\sum_{i=1}^{2n} J(\nabla_{ \nabla_{e_i} e_i} J)(X).
						\end{equation}   
						Comparing \eqref{eq:wood-1} with \eqref{eq:wood-2} we obtain 
						$[J,\nabla^* \nabla J]=2\nabla_{\delta J}J-2[J,P]$, as we wanted to prove.
					\end{proof}

					\ 
					
				\end{document}